\documentclass[11pt]{amsart} 
\usepackage{graphicx}
\usepackage[parfill]{parskip}  

\usepackage{amsthm}

\usepackage{textcomp}
\usepackage{blkarray}
\usepackage[scaled=.85]{beramono}% a typewriter font must be defined
\usepackage[leqno]{amsmath} %use eq no on left
\usepackage{bm}% load after all math to give access to bold math
\usepackage{url}
\newcommand{\Ki}{\mathcal{X}}
\newcommand{\bF}{{\mathbb{F}}}
\newcommand{\bC}{{\mathbb{C}}}
\newcommand{\bR}{{\mathbb{R}}}   
\newcommand{\bP}{{\mathbb{P}}}

\newcommand{\bL}{{\mathbb{L}}}

\newcommand{\bB}{{\mathbb{B}}}
\newcommand{\bD}{{\mathbb{D}}}
\def\coloneq{:=}
\def\Mid{\,\vert\,}
\usepackage[paperwidth=7in, paperheight=10in, textwidth=4.75in, textheight=8in, includehead=true, bindingoffset=.5in]{geometry}

\newcommand{\beq}{\begin{equation}}
\newcommand{\eeq}{\end{equation}}
\numberwithin{equation}{section}

\newtheorem{theorem}[equation]{Theorem}
  
\newtheorem{lemma}[equation]{Lemma}

\theoremstyle{remark}
\newtheorem{srem}[equation]{}

\newtheorem{remark}[equation]{Remark}

\theoremstyle{definition}
\newtheorem{defn}[equation]{Definition}
\newtheorem{definition}[equation]{Definition}

\title{\textsf{Schur Inequalities}}%\\[15pt]{\em proof techniques and intuition}}
\author{S. Gill Williamson}
\thanks{Department of Computer Science and Engineering, 
University of California San Diego; \url{http://cse.ucsd.edu/~gill}.
}
\date{}                                           % Activate to display a given date or no date
\usepackage[pdftex,colorlinks=true, pdfstartview=FitV, linkcolor=blue, citecolor=blue, urlcolor=blue,bookmarks,bookmarksnumbered]{hyperref}
\hypersetup{
pdftitle={Matrix Theory}, 
pdfauthor={S. Gill Williamson}%{\textcopyright\ S. Gill Williamson}
}
\begin{document}
\thispagestyle{empty}
\begin{center}
\vspace*{1in}
\textsf{ \Large Tensor spaces -- the basics}\\[.2in]
%\textsf{\small basic skills and proof techniques}\\[.7in]

\textsf{S. Gill Williamson}
%\vfill
%\textcopyright  \textsf{S. Gill Williamson 2012. All rights reserved.}
\end{center}

%\newpageI
\thispagestyle{empty}
\hspace{1 pt}
%\newpage
%ABSTRACT
\begin{center}
{\Large Abstract}\\[.2in]
\end{center}
\pagestyle{plain}
We present the basic concepts of tensor products of vectors spaces, exploiting the special properties of vector spaces as opposed to more general modules.
Introduction~(\ref{sec:introduction}),  
Basic multilinear algebra~(\ref{sec:basic}), 
Tensor products of vector spaces~(\ref{sec:tensorproducts}),
Tensor products of matrices~(\ref{sec:tensormatrices}),
Inner products on tensor spaces~(\ref{sec:innerproducts}),
Direct sums and tensor products~(\ref{sec:directsums}),
Background concepts and notation~(\ref{sec:review}),
Discussion and acknowledgements~(\ref{sec:discussion}).
This material draws upon \cite{mm:fdma1} and relates to subsequent work
\cite{mm:fdma2}.
%INTRODUCTION
\section{Introduction}
\label{sec:introduction}
We start with an example.
Let $x=\left(\begin{smallmatrix}a\\b\end{smallmatrix}\right)\in \mathbf{M}_{2,1}$ and 
$y=\left(\begin{smallmatrix}c&d\end{smallmatrix}\right)\in \mathbf{M}_{1,2}$
where $\mathbf{M}_{m,n}$ denotes the $m \times n$ matrices over the real numbers, $\bR$.  

The set of  matrices $\mathbf{M}_{m,n}$ forms a {\em vector space} under matrix addition and multiplication by real numbers, $\bR$. The {\em dimension} of this vector space is $mn$.  We write $\dim(\mathbf{M}_{m,n})=mn.$

Define a function $\nu(x,y)=xy=
\left(
\begin{smallmatrix}
ac\,&ad\\
bc\,&bd
\end{smallmatrix}\right)
\in \mathbf{M}_{2,2}$
(matrix product of $x$ and $y$).
The function $\nu$ has {\em domain} $V_1\times V_2$ where 
$V_1=\mathbf{M}_{2,1}$ and $V_2=\mathbf{M}_{1,2}$ are vector spaces of dimension, $\dim(V_i) = 2$, $i=1,2$.
The range of $\nu$ is the vector space  $P=\mathbf{M}_{2,2}$ which has
$\dim(P)=4.$

The function $\nu$ is {\em bilinear} in the following sense:
\begin{equation}
\label{eq:blnrxmp1}
\nu(r_1x_1 + r_2x_2, y) = r_1\nu(x_1,y)+r_2\nu(x_2,y)
\end{equation}
\begin{equation}
\label{eq:blnrxmp2}
\nu(x, r_1y_1 + r_2y_2) = r_1\nu(x,y_1)+r_2\nu(x,y_2)
\end{equation}
for any $r_1,\,r_2\in \bR,\;$ $x, x_1, x_2 \in V_1,\;$ and $y, y_1, y_2 \in V_2.$
We denote the set of {\em all} such bilinear functions by $M(V_1, V_2 : P)$.
Recall that the {\em image} of $\nu$ is the set 
$\mathrm{Im}(\nu)\coloneq\{\nu(x,y)\Mid (x,y)\in V_1\times V_2\}$ and the {\em span} of the image of $\nu$, denoted by $\left<\mathrm{Im}(\nu)\right>$,
is the set of all linear combinations of vectors in $\mathrm{Im}(\nu)$.

Let $E_1= \{e_{11}, e_{12}\}$ be an ordered basis for $V_1$ and  
$E_2 = \{e_{21}, e_{22}\}$ be an ordered basis for $V_2$
specified as follows: 
\begin{equation}
\label{eq:bsselmv1v2}
e_{11}=\left(\begin{smallmatrix}1\\0\end{smallmatrix}\right),\;
e_{12}=\left(\begin{smallmatrix}0\\1\end{smallmatrix}\right),\;
e_{21}=\left(\begin{smallmatrix}1&0\end{smallmatrix}\right),\;
e_{22}=\left(\begin{smallmatrix}0&1\end{smallmatrix}\right).
\end{equation}
Using matrix multiplication, we compute
\begin{equation}
\label{eq:nuonbasispairs}
\begin{matrix}
\nu(e_{11},e_{21}) = e_{11}e_{21} =
\left(
\begin{smallmatrix}
1\,&0\\
0\,&0
\end{smallmatrix}
\right),\;\;
\nu(e_{11},e_{22}) = e_{11}e_{22} =
\left(
\begin{smallmatrix}
0\,&1\\
0\,&0
\end{smallmatrix}
\right),\;
\\
\nu(e_{12},e_{21}) = e_{12}e_{21} =
\left(
\begin{smallmatrix}
0\,&0\\
1\,&0
\end{smallmatrix}
\right),\;\;
\nu(e_{12},e_{22}) = e_{12}e_{22} =
\left(
\begin{smallmatrix}
0\,&0\\
0\,&1
\end{smallmatrix}
\right).
\end{matrix}
\end{equation}
Note that the set 
$
\{\nu(e_{11},e_{21}), \nu(e_{11},e_{22}), \nu(e_{12},e_{21}), \nu(e_{12},e_{22})\}
\subset\mathrm{Im}(\nu)
$
is a basis for $P$.
Since $\mathrm{Im}(\nu)$ contains a basis for $P$, we have 
$\left<\mathrm{Im(\nu)}\right>=P$ (i.e., the span of the image of  $\nu$ equals $P$).
Here is our basic definition of the tensor product of two vector spaces:
%DEF
\begin{definition}[\bfseries The tensor product of two vector spaces]
\label{def:tnsprdtwovctspcs}
Let $V_1$ and $V_2$ be vector spaces over 
$\bR$ with $\dim(V_1) = n_1$ and  $\dim(V_2) = n_2$.   
Let $P$ be a vector space over $\bR$ and $\nu\in M(V_1, V_2 : P).$
Then, $(P,\nu)$ is a tensor product of the $V_1, V_2$ if
\begin{equation}
\label{eq:tnsrdfntnAsmp}
\;\;\left<\mathrm{Im}(\nu)\right> = P\;\; \mathrm{and}
\;\;\dim(P)=\prod_{i=1}^2\dim(V_i).\;\;
\end{equation}
We sometimes use
$x_1\otimes x_2\coloneq \nu(x_1,x_2)$ and call  $x_1\otimes x_2$ the {\em tensor product} of $x_1$ and $x_2$.  We sometimes use $V_1\otimes V_2$ for $P$.
\end{definition}
%SREM
\begin{srem}[\bfseries Remarks and intuition about tensor products]
Note that our  example,   $\nu(x,y)=xy=
\left(
\begin{smallmatrix}
ac\,&ad\\
bc\,&bd
\end{smallmatrix}\right)
\in \mathbf{M}_{2,2}$ 
defines a tensor product $(\mathbf{M}_{2,2}, \nu)$ by definition~\ref{def:tnsprdtwovctspcs}.
Using ``tensor'' notation, we could write 
$
\left(\begin{smallmatrix}a\\b\end{smallmatrix}\right) \otimes 
(\begin{smallmatrix}c&d\end{smallmatrix}) =
\left(
\begin{smallmatrix}
ac\,&ad\\
bc\,&bd
\end{smallmatrix}
\right)
$
and
$\mathbf{M}_{2,1}\otimes \mathbf{M}_{1,2}=\mathbf{M}_{2,2}$.
A tensor product is simply a special bilinear (multilinear) function $\nu$.  
To specify any  function you must specify its domain and range.  
Thus the $V_1$, $V_2$, $\mathrm{Im}(\nu)$ and, hence, $\left<\mathrm{Im}(\nu)\right>$ are defined when you specify $\nu$. 
In this sense the definition of a tensor needs only $\nu$, not the pair $(P,\nu)$.
It is useful, however, to display $P$ explicitly. The condition 
$\langle\mathrm{Im}(\nu)\rangle = P$ states that $P$ isn't too big in some sense---every bit of $P$ has some connection to $\nu$.  
The condition $\dim(P)=\prod_{i=1}^2\dim(V_i)$ states that $P$ isn't too small in some sense.
The function $\nu$ has enough room to reveal all of its tricks.
The notation $ x_1\otimes x_2$ is a convenient way to work with $\nu$ and not have to keep writing ``$\nu$'' over and over again.  
The choice of the {\em model} $P$ is important in different applications.
The subject of multilinear algebra (or tensors) contains linear algebra as a special case, so it is not a trivial subject.
\end{srem}

%SECTION
\section{Basic multilinear algebra}
\label{sec:basic}
Section~\ref{sec:review} gives background material for review as needed.
We focus on multilinear algebra over finite dimensional vector spaces.  
The default field, $\bF$, has characteristic zero although much of the material is valid for finite fields.
Having committed to this framework, we make use of the fact that vector spaces have bases.
Tensors are special types of multilinear functions so we need to get acquainted with multilinear functions.

%\begin{srem}[\bfseries Familiar fields]
%\label{srem:familiarfields}
%Some familiar fields, $\bF$, include $\bQ$ (rational numbers), $\bR$ (real numbers), 
%$\bC$ (complex numbers) and the fields of rational functions (ratios of polynomials) 
%$\bQ(z)$, $\bR(z)$ and $\bC(z)$.  
%Thus,  $\bF\in\{\bQ, \bR, \bC, \bQ(z), \bR(z),\bC(z)\}$. 
%\end{srem}

\begin{defn}[\bfseries  Summary of vector space and algebra axioms]
\label{def:vspaceaxioms} 
Let $\bF$ be a field and let
$(M, + )$ be an abelian group with identity $\theta$.  
Assume there is an operation,  $\bF\times M\rightarrow M$, which takes $(r,x)$ to $rx$ (juxtaposition of $r$ and $x$).  
To show that $(M, + )$ is a {\em vector space over} $\bF$, we show the following four things
hold for every $r,s \in \bF$ and $x,y \in M$:
$${\bf(1)}\;r(x+ y) = rx+  ry\;\;{\bf(2)}\;(r+s)x = rx+  sx\;\;{\bf(3)}\;(rs)x=r(sx)\;\;{\bf(4)}\;1_\bF\,x=x$$
where $1_\bF$ is the multiplicative identity in $\bF$.  If $(M, + , \cdot )$ is a ring for which 
$(M, + )$ is a vector space over $\bF$, then $(M, + , \cdot )$ is an {\em algebra over} $\,\bF$ if the following scalar rule holds: ${\bf(5)}$  for all $\alpha \in \bF$, $x, y \in M$,  
$\alpha(xy)=(\alpha x)y = x(\alpha y)$.\\
\end{defn}

\begin{srem}[\bfseries The free $\bF$-vector space over a set $K$]
\label{srem:frevctspc}
We use the notation $A^B$ to denote all functions with domain $B$ and range $A$.
Let $\bF$ be a field and $K$ a finite set.  
For  $f, g\in\bF^K$, $x\in K$, define $f+g$ by 
$(f+g)(x)\coloneq f(x)+g(x)$.
For $r \in \bF$, $f \in \bF^K$, $x\in K$, define  $(r f)(x)\coloneq r f(x)$.
It is easy to check that $(F^K,+)$ is a vector space over $\bF$ (see~\ref{def:vspaceaxioms}). 
The set of {\em indicator functions} $\{\iota_k\Mid k\in K\}$ where $\iota_k(x) = 1$ if $x=k$, $0$ otherwise, forms a basis for $\bF^K$. 
Thus, $\bF^K$ is a vector space of dimension $|K|$ (the cardinality of $K$).
If a set  $K$ has some special interest for us (e.g, elements of a finite group, finite sets of graphs or other combinatorial objects) then studying this ``free'' $F$-vector space over $K$ sometimes yields new insights about $K$.

\end{srem}

\begin{srem}[\bfseries Frequently used notations]
\label{srem:tnsprdsts}
Let $\Ki(statement) = 1$ if $statement$ is true, $0$ if false.
For short, $\delta_{\alpha, \gamma}\coloneq\Ki(\alpha = \gamma).$ 
Let $E_i=\{e_{i1}, e_{i2}, \ldots, e_{i\,n_i}\}$, $i=1,\ldots, m$, be linearly ordered sets. The order is specified by the second indices, $1, 2, \ldots, n_i$, thus, alternatively $E_i = (e_{i1}, e_{i2}, \ldots, e_{i\,n_i})$.
For any integer $m$, we define ``underline'' notation
by $\underline{m}\coloneq \{1,\ldots, m\}$.
Let 
\[\Gamma(n_1, \ldots, n_m)\coloneq 
\{\gamma\Mid
\gamma = (\gamma(1), \ldots, \gamma(m)), 1\leq \gamma(i)\leq n_i, i=1, \ldots, m
\}.
\]

The (Cartesian) product of the sets $E_i$ is 
\begin{equation}
\label{eq:tnsprdsts}
E_1\times \cdots\times E_m \coloneq \{(e_{1\gamma(1)}, \ldots, e_{m\gamma(m)})\Mid \gamma\in \Gamma(n_1, \ldots, n_m)\}
\end{equation}
where 
$\Gamma(n_1, \ldots, n_m) \coloneq \underline{n_1}\times\cdots\times \underline{n_m}.$
We use the notations $\Gamma(n_1, \ldots, n_m)$ and 
$\underline{n_1}\times\cdots\times \underline{n_m}$ as convenience dictates.
Using this notation, we have
\begin{equation}
\prod_{i=1}^2\left(\sum_{k=1}^2a_{ik}\right) = 
\sum_{\gamma\in{\underline 2}\times{\underline 2}} \prod_{i=1}^2 a_{i\,\gamma(i)}
=\sum_{\gamma\in\Gamma(2,2)} \prod_{i=1}^2 a_{i\,\gamma(i)}.
\end{equation}
The general form of this identity is
\begin{equation}
\label{eq:changesumprod}
\prod_{i=1}^m\left(\sum_{k=1}^{n_i} a_{ik}\right) = 
\sum_{\gamma\in\Gamma(n_1, \ldots, n_m)} \prod_{i=1}^m a_{i\,\gamma(i)}.
\end{equation}
This product-sum-interchange, identity \ref{eq:changesumprod}, will be used frequently in what follows.
\end{srem}

%DEF
\begin{definition}
\label{def:mltlnrfnc}
Let $V_i$, $i=1, \ldots, m$, and $U$ be vector spaces over a field $\bF$.
Let $\times_1^m V_i$ be the Cartesian product of the $V_i$.
A function $\varphi \colon \times_1^m V_i \rightarrow U$ is {\em multilinear}
if for all $t\in\underline{m}.$ 
\begin{equation*}
\varphi(\ldots, cx_t + dy_t, \ldots) = c\varphi(\ldots, x_t, \ldots)  + d\varphi( \ldots, y_t, \ldots)
\end{equation*}
where $x_i, y_i \in V_i$  and $c, d\in \bF$.
The set of  all such multilinear functions is denoted by $M(V_1, \ldots, V_m:U)$.
If $m=1$, $M(V_1:U)\coloneq \bL(V_1, U)$, the set of linear functions from $V_1$ to $U$.
\end{definition}
%SREM

\begin{srem}[\bfseries Multilinear expansion formula]
\label{srem:tnrprdint}
Let $E_i = \{e_{i1},\ldots, e_{in_i}\}$ be an ordered basis for $V_i$, $i=1,\ldots, m$.
If $\varphi\in M(V_1, \ldots, V_m:U)$ and $x_i = \sum_{k=1}^{n_i} c_{ik} e_{ik}$, $i=1,\ldots m$,
we can, using the rules of multilinearity, compute $\varphi(x_1,\dots,x_m) = $
\begin{equation}
\label{eq:mltlnintro}
\varphi\left(\sum_{k=1}^{n_1} c_{1k} e_{1k},\,\dots,\sum_{k=1}^{n_m} c_{mk} e_{mk}\right) = \sum_{\gamma\in\Gamma} \left(\prod_{i=1}^{m} c_{i\gamma(i)}\right)
\varphi(e_{1\gamma(1)},\ldots, e_{m\gamma(m)})
\end{equation}
where $\Gamma\coloneq\Gamma(n_1, \ldots, n_m)$ (see~\ref{srem:tnsprdsts}).
Formula~\ref{eq:mltlnintro} can be proved by induction. It is a generalization of the standard algebraic rules for interchanging sums and products (\ref{eq:changesumprod}).
\end{srem}
%SREM
A function with domain $\Gamma(n_1, \ldots, n_m)$ and range $S$ is sometimes denoted by ``maps to'' notation: 
 $\gamma\mapsto f(\gamma)\,$  or   $\gamma\mapsto s_{\gamma}$. 
 In the former case, $f$ is the name of the function so we write $\mathrm{Im}(f)$ for the set
 $\{f(\gamma)\Mid \gamma\in \Gamma\}$.  In the latter case (index notation), the $s$ is not consistently  interpreted as the name of the function.
 %DEF
\begin{definition}[\bfseries Multilinear function defined by extension of $\gamma\mapsto s_{\gamma}$]
\label{def:dftnsprdintro}
Let $V_i$, $i=1, \ldots, m$, be finite dimensional vector spaces over a field $F$.
Let $E_i = \{e_{i1},\ldots, e_{in_i}\}$ be an ordered basis for $V_i$, $i=1,\ldots, m$.
Let $U$ be a vector space over $\bF$  and let 
$\gamma\mapsto s_{\gamma}$ be any function with domain 
$\Gamma(n_1, \ldots, n_m)$ and range $U$.
If $x_i = \sum_{k=1}^{n_i} c_{ik} e_{ik}$, $i=1,\ldots m$, then we can use~\ref{eq:mltlnintro} to define $\nu\in M(V_1, \ldots, V_m:U)$ by
\begin{equation}
\label{eq:mltlnintro2}
\nu(x_1,\dots,x_m) \coloneq 
\sum_{\gamma\in\Gamma} \left(\prod_{i=1}^{m} c_{i\gamma(i)}\right)s_{\gamma}.
\end{equation}
We say 
$\nu\in M(V_1, \ldots, V_m:U)$ is
{\em defined by multilinear extension from $\gamma\mapsto s_{\gamma}$
and the bases $E_i$, $i=1,\ldots, m$}.

%\ref{def:tnsrprdfld}.
\end{definition}
%\SREM
\begin{srem}[\bfseries Examples of multilinear functions defined by bases]
\label{srem:exmmltfncdefbss}
We follow definition~\ref{def:dftnsprdintro}.
Let $V_1=\bR^2$ and $V_2=\bR^3$, $\bR$ the real numbers.
Let $E_1=\{e_{11}, e_{12}\}$ be the ordered basis for $\bR^2$, $e_{11}=(1,0)$, $e_{12}=(0,1).$
Let $E_2=\{e_{21}, e_{22}, e_{23}\}$ be the ordered basis for $\bR^3$, 
$e_{21}=(1,0,0)$, $e_{22}=(0,1, 0)$, $e_{23}=(0,0,1)$.
From equation~\ref{eq:tnsprdsts}, we have 
$
E_1\times E_2 =\{(e_{1\gamma(1)}, e_{2\gamma(2)})\Mid \gamma\in \Gamma(2,3)\}.
$
Let $U = \mathbf{M}_{2,3}(\bR)$.
Define an injection $\gamma\mapsto b_{\gamma} \in U$ as follows:
\begin{equation}
\label{eq:ordbasB1}
\left\{
\begin{pmatrix}
1&0&0\\
0&0&0
\end{pmatrix}_{11},
\begin{pmatrix}
0&1&0\\
0&0&0
\end{pmatrix}_{12},
\begin{pmatrix}
0&0&1\\
0&0&0
\end{pmatrix}_{13},
\begin{pmatrix}
0&0&0\\
1&0&0
\end{pmatrix}_{21},
\begin{pmatrix}
0&0&0\\
0&1&0
\end{pmatrix}_{22},
\begin{pmatrix}
0&0&0\\
0&0&1
\end{pmatrix}_{23}
\right\}.
\end{equation}
%ORDERED BASIS
The elements of the set~\ref{eq:ordbasB1} correspond to ordered pairs 
$(\gamma, b_{\gamma})$ in the form $(b_{\gamma})_{\gamma}$ and thus define $\gamma\mapsto b_{\gamma}$.  Using lexicographic order (\ref{def:lxcrdr}) on $\Gamma(2,3)$, 
the set~\ref{eq:ordbasB1} defines an ordered basis for $U$:
\begin{equation}
 B=
\label{eq:ordbasB1a}
\left(
\begin{pmatrix}
1&0&0\\
0&0&0
\end{pmatrix},\;
\begin{pmatrix}
0&1&0\\
0&0&0
\end{pmatrix},\;
\begin{pmatrix}
0&0&1\\
0&0&0
\end{pmatrix},\;
\begin{pmatrix}
0&0&0\\
1&0&0
\end{pmatrix},\;
\begin{pmatrix}
0&0&0\\
0&1&0
\end{pmatrix},\;
\begin{pmatrix}
0&0&0\\
0&0&1
\end{pmatrix}
\right).
\end{equation}
Write $B=(b_{11}, b_{12}, b_{13}, b_{21}, b_{22}, b_{23})$ where the $b_{ij}$ are the corresponding entries in the list $B$ of~\ref{eq:ordbasB1a}. 
Take $x_1= c_{11}e_{11} + c_{12}e_{12} =(c_{11}, c_{12}) \in \bR^2$.
Take $x_2= c_{21}e_{21} + c_{22}e_{22} + c_{23}e_{23} = (c_{21}, c_{22}, c_{23}) \in \bR^3$.
By equation~\ref{eq:mltlnintro2} we get $\nu_B(x_1,x_2)=$
\begin{equation}
\label{eq:mltlnintro2a}
c_{11}c_{21}b_{11}+c_{11}c_{22}b_{12}+c_{11}c_{23}b_{13}+
c_{12}c_{21}b_{21}+c_{12}c_{22}b_{22}+c_{12}c_{23}b_{23}. 
\end{equation}
Using the elements of $B$ from equation~\ref{eq:ordbasB1a}, the sum
shown in~\ref{eq:mltlnintro2a} becomes
\[
\nu_B(x_1,x_2)=
\begin{pmatrix}
c_{11}c_{21}\;&c_{11}c_{22}\;&c_{11}c_{23}\\
c_{12}c_{21}\;&c_{12}c_{22}\;&c_{12}c_{23}
\end{pmatrix}.
\]
The nice pattern here where the first row is $c_{11}(c_{21}, c_{22}, c_{23})$ and
the second row is $c_{12}(c_{21}, c_{22}, c_{23})$  is a consequence of the choice of correspondence between the elements of the basis $B$ and $\Gamma$.
For example, instead of equation~\ref{eq:ordbasB1} use 
\begin{equation}
\label{eq:ordbasB2}
B=
\left\{
\begin{pmatrix}
1&0&0\\
0&0&0
\end{pmatrix}_{11},
\begin{pmatrix}
0&1&0\\
0&0&0
\end{pmatrix}_{12},
\begin{pmatrix}
0&0&0\\
1&0&0
\end{pmatrix}_{13},
\begin{pmatrix}
0&0&1\\
0&0&0
\end{pmatrix}_{21},
\begin{pmatrix}
0&0&0\\
0&1&0
\end{pmatrix}_{22},
\begin{pmatrix}
0&0&0\\
0&0&1
\end{pmatrix}_{23}
\right\}.
\end{equation}
We then have the less attractive matrix representation
\[
\nu_B(x_1,x_2)=
\begin{pmatrix}
c_{11}c_{21}\;&c_{11}c_{22}\;&c_{12}c_{21}\\
c_{11}c_{23}\;&c_{12}c_{22}\;&c_{12}c_{23}
\end{pmatrix}.
\]
The choice of indexing (i.e., bijective correspondence to $\Gamma(2,3)$) matters.
\end{srem}
%SREM
\begin{srem}[\bfseries Component spaces of $M(V_1, \ldots, V_m:U)$ are  $M(V_1, \ldots, V_m:\bF)$]
\label{rem:vcspcfld}
Suppose that  $U=\left<u_1, \ldots, u_n\right>$, where angle brackets, ``$\left<\right>$''  denote the span of the basis $\{u_1, \ldots, u_n\}$ for $U$. 
Note that $U$ is the direct sum $\oplus^n\left<u_i\right>$ of
the one-dimensional subspaces $\left<u_i\right>$.  
The vector space $U$ is isomorphic to the vector space $\bF^{\underline{n}}$ and the $\left<u_i\right>$ are isomorphic (as vector spaces over $\bF$) to $\bF$.
The function $\varphi\in M(V_1, \ldots, V_m:U)$ has component functions 
$\varphi^{(i)}$:
\begin{equation}
\label{eq:cmpfncs}
\varphi(x_1, \ldots, x_m) = 
\sum_{i=1}^n\varphi^{(i)}(x_1, \ldots, x_m)u_i.
\end{equation}
Each $\varphi^{(i)} \in M(V_1, \ldots, V_m:\left<u_i\right>)\equiv M(V_1, \ldots, V_m:\bF)$.
\end{srem}
%REM
\begin{srem}[\bfseries Vector space of multilinear functions]
\label{rem:mltlnvcspc}
The set of all multilinear functions, $M(V_1, \ldots, V_m:U)$, is a vector space over $\bF$ in the standard way: 
\[(\varphi_1 + \varphi_2)(x_1, \ldots, x_m)=\varphi_1(x_1, \ldots, x_m)+
\varphi_2(x_1, \ldots, x_m)
\]
and $(c\varphi)(x_1, \ldots, x_m) = c\varphi(x_1, \ldots, x_m).$
Suppose, for $i\in \underline{m}$, $V_i = \left<E_i\right>$ ,  $E_i = \{e_{i1},\ldots, e_{in_i}\}$ a 
basis for $V_i$.
Let $U=\left<u_1, \ldots, u_n\right>$ where $\{u_i\Mid i\in \underline {n}\}$ is a
basis.
Define the component functions~(\ref{eq:cmpfncs}) $\varphi^{(i)} \in 
M(V_1, \ldots, V_m:\left<u_i\right>)$ by 
\begin{equation}
\label{eq:bssmltlnfncs}
\varphi^{(j)}_{\alpha}(e_{1\gamma(1)},\ldots, e_{m\gamma(m)}) \coloneq 
\delta_{\alpha, \gamma}u_j.
\end{equation}
The set of $n\prod_1^m n_i$ multilinear functions 
\begin{equation}
\label{eq:mltlnrfncbss00}
\{\varphi^{(j)}_{\alpha}\Mid \alpha \in \Gamma(n_1,\ldots, n_m),\; j\in \underline{n}\}
\end{equation}
is  a basis for the vector space $M(V_1, \ldots, V_m:U).$ 
See~\ref{rem:mltlnrfncbss} for discussion.
Note that as vector spaces
\begin{equation}
\label{eq:multdirsum1}
M(V_1, \ldots, V_m:U)=\bigoplus_{j=1}^n M(V_1, \ldots, V_m:\left<u_j\right>).
\end{equation}
\end{srem}

%REM
\begin{srem}[\bfseries The multilinear function basis]
\label{rem:mltlnrfncbss}
It suffices to take $U=\bF$ (see~\ref{rem:vcspcfld}).
We show
$\{\varphi_{\alpha}\Mid \alpha \in \Gamma(n_1,\ldots, n_m)\}$
is a basis for $M(V_1, \ldots, V_m:\bF)$ where
$
\varphi_{\alpha}(e_{1\gamma(1)},\ldots, e_{m\gamma(m)}) =  \delta_{\alpha, \gamma}1_{\bF}.
$

Here, $E_i = \{e_{i1},\ldots, e_{in_i}\}$ as in~\ref{rem:mltlnvcspc}.
If $\varphi\in M(V_1, \ldots, V_m:\bF)$ and $x_i = \sum_{k=1}^{n_i} c_{ik} e_{ik}$
we have
\begin{equation}
\label{eq:mltlnrfncbss01}
\varphi(x_1,\dots,x_m) = \sum_{\gamma\in\Gamma} \left(\prod_{i=1}^{m} c_{i\gamma(i)}\right)
\varphi(e_{1\gamma(1)},\ldots, e_{m\gamma(m)}).
\end{equation}
Replacing $\varphi$ by $\varphi_{\alpha}$ in equation~\ref{eq:mltlnrfncbss01} we obtain   
\begin{equation}
\label{eq:mltlnrfncbss1}
\varphi_{\alpha}(x_1,\ldots, x_m) = \sum_{\gamma\in\Gamma} 
\left(\prod_{i=1}^{n_i} c_{i\gamma(i)}\right) 
\varphi_{\alpha}(e_{1\gamma(1)},\ldots, e_{m\gamma(m)})=\prod_{i=1}^{m} c_{i\alpha(i)}.
\end{equation}
From equations~\ref{eq:mltlnrfncbss01} and~\ref{eq:mltlnrfncbss1} we obtain 
\begin{equation}
\label{eq:mltlnrfncbss2}
\varphi(x_1,\dots,x_m) = \sum_{\alpha\in\Gamma} 
\varphi_{\alpha}(x_1,\ldots, x_m)
\varphi(e_{1\alpha(1)},\ldots, e_{m\alpha(m)}).
\end{equation}
Thus, from \ref{eq:mltlnrfncbss2}, we obtain
\begin{equation}
\label{eq:mltlnrfncbss3}
\varphi = \sum_{\alpha\in\Gamma} \varphi(e_{1\alpha(1)},\ldots, e_{m\alpha(m)}) \varphi_{\alpha}
\end{equation}
which can be proved by evaluating both sides at $(x_1, \ldots, x_m)$.
It is obvious that $\{\varphi_{\alpha}\Mid \alpha \in \Gamma(n_1,\ldots, n_m)\}$
is linearly independent for if 
$\varphi=\sum_{\alpha} d_{\alpha} \varphi_{\alpha}=0$ then 
$\varphi(e_{1\gamma(1)},\ldots, e_{m\gamma(m)}) = d_{\gamma} = 0$ for all $\gamma$.
\end{srem}
%
%SECTION
\section{Tensor products of vector spaces}
\label{sec:tensorproducts}
%DEF
\begin{definition}[\bfseries Tensor product]
\label{def:tnsrprdfld}
Let $V_1, \ldots, V_m$, be vector spaces over $\bF$ of dimension $\dim(V_i)$, $i\in \underline{m}$.
Let $P$ be a vector space over $\bF$. 
Then the pair $(P, \nu)$ is a {\em tensor product} of the $V_i$ if
\begin{equation}
\label{eq:tnsrdfntn0}
\nu\in M(V_1, \ldots, V_m:P)\;\; \mathrm{and}
\;\;\left<\mathrm{Im}(\nu)\right> = P\;\; \mathrm{and}
\;\;\dim(P)=\prod_{i=1}^m\dim(V_i)
\end{equation}
where $\left<\mathrm{Im}(\nu)\right>\coloneq
\left<\nu(x_1, \ldots, x_m)\Mid (x_1, \ldots, x_m)\in \times_{i=1}^m V_i\right>$
is the span of $\mathrm{Im}(\nu)$, the {\em image} of $\nu$.
The vector space $P$ is denoted by $V_1\otimes\cdots\otimes V_m$. The vectors
$\nu(x_1,\ldots, x_m)\in P$ are designated by $x_1\otimes \cdots \otimes x_m$ and are called the {\em homogeneous} tensors.
The notation $\otimes(x_1, \ldots , x_m) \coloneq \nu(x_1, \ldots , x_m)$ is sometimes used.
\end{definition}
%SREM
%LEM
\begin{lemma}[\bfseries Canonical bases]
\label{lem:cnnbss}
Let $V_1, \ldots, V_m$, $\dim(V_i) = n_i$, $i\in \underline{m}$, be vector spaces over $\,\bF$. 
Let $P$ be a vector space over $\,\bF$ and $\nu\in M(V_1, \ldots, V_m:P).$
Suppose, for all $i\in \underline{m}$, $V_i = \left<E_i\right>$ ,  
$E_i = \{e_{i1},\ldots, e_{in_i}\}$ an ordered basis for $V_i$.
Then, $(P,\nu)$ is a tensor product of the $V_i$ if and only if 
\begin{equation}
\label{eq:cnnbss}
\bB\coloneq\{\nu(e_{1\gamma(1)},\ldots, e_{m\gamma(m)})\Mid \gamma\in 
\underline{n_1}\times \cdots\times \underline{n_m}\} \;\mathit{is\;a\;basis\;for}\; P.
\end{equation}
The standard convention is that $\bB$ is ordered using lexicographic order
on $\Gamma$. 
\begin{proof}
Suppose $(P,\nu)$ is a tensor product and assume $x_i = \sum_{k=1}^{n_i} c_{ik} e_{ik}$ so that
\begin{equation}
\label{eq:cnnbss2}
\nu(x_1,\dots,x_m) = \sum_{\gamma\in\Gamma} \left(\prod_{i=1}^{m} c_{i\gamma(i)}\right)
\nu(e_{1\gamma(1)},\ldots, e_{m\gamma(m)}).
\end{equation}
By definition~\ref{def:tnsrprdfld},  $\left<\mathrm{Im}(\nu)\right>=
\left<\nu(x_1, \ldots, x_m)\Mid (x_1, \ldots, x_m)\in \times_{i=1}^m V_i\right>=P.$
Thus, $\bB =
\{\nu(e_{1\gamma(1)},\ldots, e_{m\gamma(m)})\Mid \gamma\in 
\underline{n_1}\times \cdots\times \underline{n_m}\}\;\mathrm{spans}\;P$.
Thus, $|\bB|=\dim(P)$ and  $\bB$ is a basis for $P$.
Conversely, if $\bB$ is a basis for $P$ then clearly $\left<\mathrm{Im}(\nu)\right>=P$ and 
$\dim(P)=|\bB|=\prod_{i=1}^mn_i$.
\end{proof}
\end{lemma}
%REM
\begin{srem}[\bfseries The space $P$ of $(P,\nu)$  
can be any $P$ with $\dim(P)=\prod_i \dim(V_i)$]
\label{rem:tnsprdalwexs} 
For short, we write  $\Gamma=\Gamma(n_1,\ldots, n_m).$
As before, $\mathrm{dim}(V_i)=n_i$. 
Let $P$ be any vector space of dimension $\prod_{i=1}^mn_i$. Let 
$\{p_{\gamma} \Mid \gamma\in \Gamma\}$ be any basis for $P$.  
Let $E_i = \{e_{i1},\ldots, e_{in_i}\}$ be the bases for $V_i$, $i=1,\ldots, m$.
Define 
$\nu(e_{1\gamma(1)},\ldots, e_{m\gamma(m)})=p_{\gamma}$ for $\gamma\in\Gamma$
and  define $\nu\in M(V_1, \ldots, V_m:P)$ by multilinear extension:
\begin{equation}
\label{eq:mltlnrfncb}
\nu(x_1,\dots,x_m) = \sum_{\gamma\in\Gamma} \left(\prod_{i=1}^{m} c_{i\gamma(i)}\right)
\nu(e_{1\gamma(1)},\ldots, e_{m\gamma(m)}).
\end{equation}
By lemma~\ref{lem:cnnbss}
the pair $(P,\nu)$ is a tensor product, $V_1, \ldots, V_m$.
The homogenous tensors are, by definition, $\nu(x_1,\dots,x_m)=
x_1\otimes\cdots\otimes x_m$ and
\[
\nu(e_{1\gamma(1)},\ldots, e_{m\gamma(m)})=
e_{1\gamma(1)}\otimes\cdots\otimes e_{m\gamma(m)}=p_{\gamma}.
\] 
In terms of homogeneous tensors, equation~\ref{eq:mltlnrfncb} becomes
\begin{equation}
\label{eq:dfnhmgtns}
x_1\otimes\cdots\otimes x_m = \sum_{\gamma\in\Gamma} \left(\prod_{i=1}^{m} c_{i\gamma(i)}\right)
e_{1\gamma(1)}\otimes\cdots\otimes e_{m\gamma(m)}.
\end{equation}

\end{srem}
 %DEF
\begin{lemma}[\bfseries Universal factorization (UF) property]
\label{lem:tnsrprdvctrsp}
Let $P$ and $V_1, \ldots, V_m$, $\dim(V_i) = n_i$, $i=1,\ldots, m$, be vector spaces over $\bF$. 
The pair $(P, \nu)$ is a {\em tensor product} of the $V_i$ if and only if
\begin{equation}
\label{eq:tnsrdfntn1}
\nu\in M(V_1, \ldots, V_m:P)\;\; \mathrm{and}
\;\;\left<\mathrm{Im}(\nu)\right> = P\;\; \mathrm{and}
\end{equation}
$\mathbf{UF}$: For any  $\varphi\in M(V_1, \ldots, V_m:\bF)$ there exists 
$h_{\varphi}\in \bL(P,\bF)$ with $\varphi=h_{\varphi}\nu.$
The statement $\mathbf{UF}$ is called the {\em Universal Factorization} property.
\begin{proof}
First assume $(P, \nu)$ is a {\em tensor product} of the $V_i$.
The set 
\[
\bB=\{\nu(e_{1\gamma(1)},\ldots, e_{m\gamma(m)})\Mid \gamma\in \Gamma\}
\]
is a basis for $P$ by lemma~\ref{lem:cnnbss}. 
For $\varphi\in M(V_1, \ldots, V_m:F)$, define   
\begin{equation}
\label{eq:mltext2}
h_{\varphi}\left(\nu(e_{1\gamma(1)},\ldots, e_{m\gamma(m)})\right) = 
\varphi(e_{1\gamma(1)},\ldots, e_{m\gamma(m)}).
\end{equation} 
Thus,  $\varphi$ and $h_{\varphi}\nu\in M(V_1, \ldots, V_m:\bF)$  agree on 
the $(e_{1\gamma(1)},\ldots, e_{m\gamma(m)})\in \times_{i=1}^mV_i$
and, by multilinear extension, $h_{\varphi}\nu = \varphi$.
The function $h_{\varphi}$ is defined on the basis $\bB$ and thus is in 
$\bL(P,\bF)$ by linear extension.

Assume conditions~\ref{eq:tnsrdfntn1} and $\mathbf{UF}$ hold. 
Let $\{\varphi_{\alpha}\Mid \alpha \in \Gamma(n_1,\ldots, n_m)\}$
be the standard basis for $M(V_1, \ldots, V_m:\bF)$ (\ref{rem:mltlnrfncbss}). 
By assumption, for each $\varphi_{\alpha}$ there exists $h_{\varphi_{\alpha}}\in \bL(P,\bF)\coloneq P^*$  (dual space to $P$) such that
%$h_{\varphi_{\alpha}}\nu = \varphi_{\alpha}$.
\begin{equation}
\label{eq:relatetodual}
\delta_{\alpha, \gamma} = \varphi_{\alpha}(e_{1\gamma(1)},\ldots, e_{m\gamma(m)}) =
h_{\varphi_{\alpha}}\nu(e_{1\gamma(1)},\ldots, e_{m\gamma(m)}). 
\end{equation} 

To show that $\bB=\{\nu(e_{1\gamma(1)},\ldots, e_{m\gamma(m)})\Mid \gamma\in 
\Gamma(n_1,\ldots, n_m)\}$ is linearly independent, assume
\[
\sum_{\gamma\in \Gamma}k_{\gamma} \nu(e_{1\gamma(1)},\ldots, e_{m\gamma(m)})=0.
\]
For  $\alpha\in \Gamma$ we use~\ref{eq:relatetodual} to obtain
$h_{\varphi_{\alpha}}\left(\sum_{\gamma\in \Gamma}k_{\gamma} \nu(e_{1\gamma(1)},\ldots, e_{m\gamma(m)}\right) = $
\[
\sum_{\gamma\in \Gamma}k_{\gamma} h_{\varphi_{\alpha}}\nu(e_{1\gamma(1)},\ldots, e_{m\gamma(m)}) =\sum_{\gamma\in \Gamma}k_{\gamma}\delta_{\alpha,\gamma} =k_{\alpha}=0.
\]
 Thus,  
\begin{equation}
\label{eq:cnnbss3}
\bB=\{\nu(e_{1\gamma(1)},\ldots, e_{m\gamma(m)})\Mid \gamma\in 
\Gamma(n_1,\ldots, n_m)\}
\end{equation}
is linearly independent and, by the condition 
$\langle\mathrm{Im}(\nu)\rangle=P$, spans $P$. 
Thus, $\bB$ is a basis and $(P,\nu)$ is a tensor product by lemma~\ref{lem:cnnbss}.   
The set  $\{h_{\varphi_{\alpha}}\Mid \alpha\in\Gamma\}$ 
is  the dual basis to $\bB$.                                          
\end{proof}
\end{lemma}

\begin{remark}[\bfseries Remarks about the universal factorization property]
\label{rem:rplcubyf}
Conditions similar to those of lemma~\ref{lem:tnsrprdvctrsp} are commonly used to define tensor products on more general algebraic structures.
Our proof, together with discussion~\ref{rem:mltlnvcspc}, shows that for finite dimensional vector spaces, $\bF$ can be replaced by $U$ in 
lemma~\ref{lem:tnsrprdvctrsp} to give an equivalent  (for vector spaces) form of 
$\mathbf{UF}$
which is summarized by commutative 
diagram~\ref{eq:Fig001}.\\
\end{remark}
%COMMUTATIVE DIAGRAM
\noindent
\begin{minipage}{\textwidth}
\begin{equation}
\label{eq:Fig001}
\mathbf{Universal\;factorization\;diagram}
\end{equation}
\begin{center}
\includegraphics{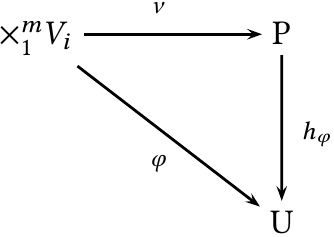}
\end{center}
\end{minipage}

\begin{remark}[\bfseries Equivalent definitions]
\label{rem:eqvdfnrtnsprd}
Let $V_1, \ldots, V_m$, $\dim(V_i) = n_i$, be vector spaces over $\bF$. 
Let $P$ be a vector space over $\bF$ and $\nu\in M(V_1, \ldots, V_m:P).$
Suppose, for all $i=1,\ldots m$,  
$E_i = \{e_{i1},\ldots, e_{in_i}\}$ is a basis for $V_i$.
Let $\Gamma = \Gamma(n_1,\ldots, n_m)$ (notation~\ref{srem:tnsprdsts}).
Then, the following statements are equivalent: 
%$(P,\nu)$ is a tensor product of the $V_i$ if and only if 
\begin{equation}
\label{eq:tnsrdfntnA}
\;\;\left<\mathrm{Im}(\nu)\right> = P\;\; \mathrm{and}
\;\;\dim(P)=\prod_{i=1}^m\dim(V_i)
\end{equation}
\begin{equation}
\label{eq:cnnbssA}
\bB=\{\nu(e_{1\gamma(1)},\ldots, e_{m\gamma(m)})\Mid \gamma\in 
\Gamma\} \;\mathit{is\;a\;basis\;for}\; P
\end{equation}
\begin{equation}
\label{eq:tnsrdfntn1A}
\;\;\left<\mathrm{Im}(\nu)\right> = P\;\; \mathrm{and}\;\;\mathbf{Universal\; Factorization:}
\end{equation}
({\bf UF}) For any  
$\varphi\in M(V_1, \ldots, V_m:\bF)$ there exists a linear function
$h_{\varphi}\in \bL(P,\bF)$ with $\varphi=h_{\varphi}\nu$.

Statement \ref{eq:tnsrdfntnA} is our definition of tensor product 
(definition~\ref{def:tnsrprdfld}).
The equivalence of statement \ref{eq:cnnbssA} and  \ref{eq:tnsrdfntnA}   follows from lemma~\ref{lem:cnnbss}.  
The equivalence of statements \ref{eq:tnsrdfntn1A} and  \ref{eq:tnsrdfntnA}     follows from lemma~\ref{lem:tnsrprdvctrsp} (Note: The proof of 
lemma~\ref{lem:tnsrprdvctrsp} shows that $h_{\varphi}$ can be assumed to be unique in this lemma).
\end{remark}
%SREM

%DEF
\begin{definition}[\bfseries Subspace tensor products]
\label{def:subtenpros}
Let $V_1, \ldots, V_m$, $\dim(V_i) = n_i$, be vector spaces.
Let $W_1\subseteq V_1, \ldots, W_m\subseteq V_m$, be subspaces. 
Let $(P, \nu)$ be a tensor product of $V_1, \ldots, V_m$.
If $(P_w, \nu_w)$ is a tensor product of $W_1, \ldots, W_m$ such that 
$\nu_w(x_1, \ldots, x_n)=\nu(x_1, \ldots, x_m)$ for 
$(x_1, \ldots, x_m)\in \times_{i=1}^m W_i$ then
$(P_w, \nu_w)$ is a {\em subspace tensor product} of $(P, \nu)$.
\end{definition}
\begin{srem}[\bfseries Subspace tensor products--constructions]
\label{srem:subtenpros}
Let $V_1, \ldots, V_m$ with dimensions $\dim(V_i) = n_i$ be vector spaces.
Let $W_1\subseteq V_1, \ldots, W_m\subseteq V_m$ be subspaces
and $(P, \nu)$ be a tensor product of $V_1, \ldots, V_m$.
Construct ordered bases $E'_i = \{e_{i1},\ldots, e_{ir_i}\}$ for the $W_i$, $i=1,\ldots, m$.
Extend these bases to $E_i = \{e_{i1},\ldots, e_{in_i}\}$, ordered bases for $V_i$, $i=1,\ldots, m$.
By \ref{eq:cnnbssA}
\begin{equation}
\label{eq:cnnbssA1}
\bB=\{\nu(e_{1\gamma(1)},\ldots, e_{m\gamma(m)})\Mid \gamma\in 
\Gamma(n_1, \ldots, n_m)\} \;\mathrm{is\;a\;basis\;for}\; P.
\end{equation}
Thus,
\begin{equation}
\label{eq:cnnbssA2}
P_w=\langle\{\nu(e_{1\gamma(1)},\ldots, e_{m\gamma(m)})\Mid \gamma\in 
\Gamma(r_1, \ldots, r_m)\}\rangle
\end{equation}
defines a subspace $P_w\subseteq P$ of dimension $\prod_{i=1}^m \dim(W_i).$ 
Note that for $\alpha\in \Gamma(r_1, \ldots, r_m),\;$
 $(e_{1\alpha(1)},\ldots, e_{m\alpha(m)})\in \times_{i=1}^m W_i,\,$ 
$\nu_w (e_{1\alpha(1)},\ldots, e_{m\alpha(m)}) \coloneq
 \nu (e_{1\alpha(1)},\ldots, e_{m\alpha(m)}).$ By multilinear extension,
$\nu_w \in M(W_1, \ldots, W_m : P_w)$.
By identity~\ref{eq:tnsrdfntnA} or~\ref{eq:cnnbssA}, $(P_w, \nu_w)$ is a tensor product of 
$W_1, \ldots, W_m$ where $P_w$ is a subspace of $P$ and $\nu_w$ is 
the restriction of $\nu$ to $\times_{i=1}^m W_i.$ 
\end{srem}

%SREM
%\begin{srem}[\bfseries Remarks on sub-tensor products]
%\label{srem:subtenprosexp}
%Definition~\ref{def:subtenpros} states, in effect, that $P_w$ is the subspace
%of $\otimes_{i=1}^mV_i$ spanned by 
%$\nu(x_1, \ldots, x_m)=x_1\otimes \cdots \otimes x_m$.
%Note that $\otimes_{i=1}^m W_i$ is a pair $(P_w, \nu_w)$ where
%$\nu_w \in M(W_1, \ldots, W_m: F)$.  However, $\nu\in  M(V_1, \ldots, V_m: F)$, and in general $\nu \neq \nu_w$ (different domains).  
%The construction of~\ref{srem:subtenpros} is always possible and represents the standard way sub-tensor products of vector spaces are constructed.
%\end{srem}
%SREM
\begin{srem}[\bfseries It's not hard to be a tensor product]
Equation~\ref{eq:relatetodual} and the related discussion shows the connection between
the $\mathbf{UF}$ property and the dual space, $P^{*}$, of $P$.
Referring to figure~\ref{eq:Fig001}, take $m=1$ and $U=\bF$.
Thus, $\times_1^m V_i = V_1$.  
Take $P=V_1$, $\nu = id$ (the identity function  in $\bL(V_1,V_1)$).  
Let $\varphi\in M(V_1:\bF)=\bL(V_1,\bF)=V_1^{*}$ (the dual space of $V_1$).
Take $h_{\varphi} = \varphi$.  
Thus, $(V_1, id)$ is a tensor product for $V_1$ (see figure~\ref{eq:Fig0012}).\\
\noindent
\begin{minipage}{\textwidth}
\begin{equation}
\label{eq:Fig0012}
%\mathbf{Tensor\;commutative\;diagram}
\end{equation}
\begin{center}
\includegraphics{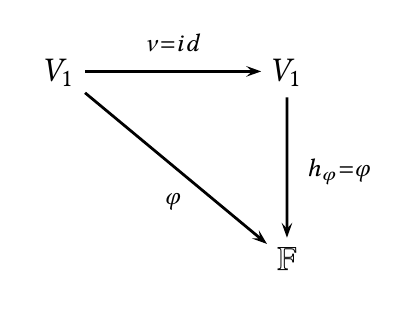}
\end{center}
\end{minipage}
\end{srem}

%REM
\begin{srem}[\bfseries The dual space model for tensors]
\label{rem:}
Referring to definition~\ref{def:tnsrprdfld}, suppose that $(P,\nu)$ is a tensor product for $V_1, \ldots, V_m$.
By~\ref{rem:tnsprdalwexs}, the vector space $P$ can be {\em any} vector space over $\bF$ with dimension 
$N=\prod_{i=1}^m\mathrm{dim}(V_i)$.  

Choose $P=M^{*}$, the dual space to the vector space
$M\coloneq M(V_1, \ldots, V_m:\bF)$.  
Define a tensor product of $V_1,\ldots, V_m$ to be $(M^{*},\nu).$  
For $\varphi\in M(V_1, \ldots, V_m:\bF)$, $\nu$ is defined by 
$\nu(x_1, \ldots, x_m)(\varphi) = \varphi(x_1, \ldots, x_m)$.
It is easy to see that $\nu$ is multilinear and well defined. 
From~\ref{rem:mltlnrfncbss}, 
$\{\varphi_{\alpha}\Mid \alpha \in \Gamma(n_1,\ldots, n_m)\}$,
where
$
\varphi_{\alpha}(e_{1\gamma(1)},\ldots, e_{m\gamma(m)}) =  \delta_{\alpha, \gamma},
$
is a basis for $M(V_1, \ldots, V_m:\bF)$.
Thus, 
\[
\{\nu(e_{1\gamma(1)},\ldots, e_{m\gamma(m)})\Mid \gamma\in\Gamma(n_1,\ldots, n_m)\}
\]
is the basis for $M^{*}$ dual to the basis $\varphi_{\alpha}$ for $M$.
The tensor product $(M^{*}, \nu)$ of $V_1, \ldots, V_n$ 
is called the
{\em dual tensor product} of $V_1,\ldots, V_m$.  As usual, 
the homogeneous tensors are $\nu(x_1, \ldots, x_m)\coloneq x_1\otimes\cdots\otimes x_m$. These homogeneous tensors satisfy  $x_1\otimes\cdots\otimes x_m(\varphi) = \varphi(x_1,\ldots, x_m)$.
\end{srem}
%DEF
\begin{definition}[\bfseries Matrix of a pair  of bases]
\label{def:matpaibas}
\index{matrix!of base pair}
\index{linear transformation!matrix given bases}
Let ${\bf v}=(v_1, \ldots v_q)$ and ${\bf w}=(w_1, \ldots, w_r)$ be ordered bases for $V$ and $W$
respectively.
Suppose for  each $j$, $1\leq j\leq q$, $T(v_j)=\sum_{i=1}^r a_{ij}w_i$. 
The matrix $A=(a_{ij})$ is called the matrix of $T$ with respect to the  base pair $({\bf v},{\bf w})$.
We write $[T]_{\bf v}^{\bf w}$ for $A$.
\end{definition}

For example, let $V = \bR^2$ and $W=\bR^3$.   
Let ${\bf w}= \{w_1, w_2, w_3\}$ and ${\bf v}=\{v_1, v_2\}$.  
Define $T$ by $T(v_1) = 2w_1 + 3w_2 - w_3$ and $T(v_2) = w_1 + 5w_2 + w_3$.
Then
\[[T]_{\bf v}^{\bf w} = \left[\begin{array}{cc} 2&1\\3&5\\-1&1\end{array}\right]\]
is the matrix of $T$ with respect to the  base pair $({\bf v},{\bf w})$.
\begin{srem}[\bfseries Obligatory observations about  tensor products of $V_1, \ldots, V_m$.]
\label{rem:alltsrprdsm}
Let $(P,\nu)$ and $(Q,\mu)$ be two tensor products of $V_1, \ldots, V_m$.
Both $P$ and $Q$ are vector spaces over $\bF$ of the same dimension,
$N=\prod_{i=1}^m \dim(V_i)$. 
Thus, they are isomorphic.
By~\ref{lem:cnnbss},
$
\bB_{\nu}\coloneq\{\nu(e_{1\gamma(1)},\ldots, e_{m\gamma(m)})\Mid \gamma\in \Gamma\}
$ 
is a basis for $P$, and
$
\bB_{\mu}\coloneq\{\mu(e_{1\gamma(1)},\ldots, e_{m\gamma(m)})\Mid \gamma\in \Gamma\}
$ 
is a  basis for $Q$.
The isomorphism $T\in \bL(P,Q)$ defined by these bases: 
\begin{equation}
\label{eq:alltsrprdsm1}
T(\nu(e_{1\gamma(1)},\ldots, e_{m\gamma(m)})\coloneq\mu(e_{1\gamma(1)},\ldots, e_{m\gamma(m)})
\end{equation} 
is a natural choice of correspondence between these tensor products.
By multilinear extension, equation~\ref{eq:alltsrprdsm1} implies 
$T\nu = \mu \in M(V_1, \ldots, V_m: Q)$.  
In matrix terms, ordering ${\bB_{\nu}}$and ${\bB_{\mu}}$ lexicographically based on $\Gamma$ gives
\[
\begin{bmatrix}
T
\end{bmatrix}_{\bB_{\nu}}^{\bB_{\mu}}
=I_N.
 \]
%Let $\varphi\in M(V_1, \ldots, V_m:F)$.
Apply lemma~\ref{lem:tnsrprdvctrsp} to $(P,\nu)$ 
with $U=Q$ and $\varphi = \mu\in M(V_1, \ldots, V_m:Q).$
Let $h_{\nu,\mu}\in \bL(P,Q)$ be such that $h_{\nu,\mu}\nu = \mu.$ 
Thus, $h_{\nu,\mu} = T$ of equation~\ref{eq:alltsrprdsm1}.
Reversing the roles of $P$, $Q$ and $\nu$, $\mu$, 
let $h_{\mu,\nu}$ be such that $h_{\mu,\nu}\mu = \nu$.
Thus, $h_{\mu,\nu} = T^{-1}\in \bL(Q,P)$.

%\begin{equation}
%\label{eq:alltsrprdsm2}
%h^{\nu}_{\mu}\nu(e_{1\gamma(1)},\ldots, e_{m\gamma(m)})=
% h^{\mu}_{\nu}\mu(e_{1\gamma(1)},\ldots, e_{m\gamma(m)})
%\end{equation} 
%so that $h^{\nu}_{\varphi}= h^{\mu}_{\varphi}T$ on the basis
%$\bB_{\nu}$.  Thus,  
%$h^{\nu}_{\varphi}= h^{\mu}_{\varphi}T= \varphi\in M(V_1, \ldots, V_m:F)$. 
\end{srem}

\begin{srem}[\bfseries Associative laws for tensor products]
\label{rem:asslwstnsprd}
Suppose, for all $i\in \underline{m}$, $V_i = \left<E_i\right>$ is the space spanned by the ordered basis $E_i = \{e_{i1},\ldots, e_{in_i}\}$.
Form a tensor product 
\[Z=(V_1\otimes\cdots\otimes V_p)\otimes(V_{p+1}\otimes\cdots\otimes V_{p+q})\]
where $p+q=m$.  
Thus, we have a tensor product of two finite dimensional vector spaces
$P_1=V_1\otimes\cdots\otimes V_p$ and $P_2=V_{p+1}\otimes\cdots\otimes V_{p+q}.$
We have $\dim(P_1)=\prod_{i=1}^p\dim(V_i)$ and $\dim(P_2)=\prod_{i=p+1}^m\dim(V_i)$.
Thus, 
\[
\dim(Z)=\dim(P_1)\dim(P_2)=\prod_1^m \dim(V_i) = \dim(V_1\otimes\cdots\otimes V_m)
\] so
$Z$ and  $(V_1\otimes\cdots\otimes V_m)$  are isomorphic vector spaces.
There are  many isomorphisms, but some are more ``natural'' than others.
Take 
$\Gamma_1(n_1, \ldots, n_p) \coloneq \underline{n_1}\times\cdots\times \underline{n_p}$ and
$\Gamma_2(n_{n_{p+1}}, \ldots, n_m)\coloneq \underline{n_{p+1}}\times\cdots\times \underline{n_m}$
and order each lexicographically.  
The natural basis for $Z$ is $e_{\alpha}\otimes e_{\beta}$ where
\[e_{\alpha}\otimes e_{\beta}=(e_{1\alpha(1)}\otimes\cdots\otimes e_{p\alpha(p)})\otimes 
(e_{(p+1)\beta(p+1)}\otimes\cdots\otimes e_{m\beta(m)})\]
where $(\alpha,\beta)\in \Gamma_1(n_1, \ldots, n_p)\times\Gamma_2(n_{n_{p+1}}, \ldots, n_m) $
is ordered lexicographically based on the lexicographic orders on $\Gamma_1$ and $\Gamma_2$.
A pair, $(\alpha,\beta)$ defines a $\gamma\in \Gamma(n_1, \ldots, n_m)$ by concatenation, 
$(\alpha,\beta)\mapsto \alpha\beta\coloneq\gamma$. This correspondence is order preserving if we assume lexicographic order on 
$\Gamma(n_1, \ldots, n_m).$
\end{srem}%
%DEF
\section{Tensor products of matrices}
\label{sec:tensormatrices}
%REM
\begin{srem}[\bfseries Tensor product of matrices example]
\label{rem:tnsprdmtr}
Let $V_1=\mathbf{M}_{p,q}$ ($p\times q$ matrices over $\bF$), and let $V_2=\mathbf{M}_{r,s}$. 
Let $A=(a_{ij})\in \mathbf{M}_{p,q}$ and $B=(b_{ij})\in \mathbf{M}_{r,s}$.
A tensor product $(P,\nu)=V_1\otimes V_2$ can be constructed following the general 
approach of~\ref{rem:tnsprdalwexs}.
Let $P$ be any vector space over $\bF$ of dimension $pqrs$.   
Let $\{p_{ijkl}\Mid (i,j,k,l)\in \underline{p}\times\underline{q}\times\underline{r}\times\underline{s}\}$
be a basis for $P$.
Let $\{E_{ij}^{(1)}\Mid (i,j)\in \underline{p}\times\underline{q}\}$ be the standard basis for $\mathbf{M}_{p,q}$ (i.e., $E_{ij}^{(1)}(i',j')=\mathcal{X}((i,j)=(i',j'))$). 
Let $\{E_{ij}^{(2)}\Mid (i,j)\in \underline{r}\times\underline{s}\}$ be the standard basis for 
$\mathbf{M}_{r,s}$.
Define $\nu( E_{ij}^{(1)},E_{kl}^{(2)})\coloneq E_{ij}^{(1)}\otimes E_{kl}^{(2)} = p_{ijkl}$, and by multilinear (bilinear here) extension
\begin{equation}
\label{eq:tnsprdtwomtr1}
\nu(A,B)\coloneq A\otimes B = \sum_{(i,j,k,l)} a_{ij}b_{kl}E_{ij}^{(1)}\otimes E_{kl}^{(2)}
 \end{equation}
 where
\begin{equation}
\label{eq:tnsprdtwomtr2}
 A = \left( \sum_{(i,j)\in \underline{p}\times\underline{q}}a_{ij}E_{ij}^{(1)}\right)\;\;
 \mathrm{and}\;\;
 B=\left( \sum_{(i,j)\in \underline{r}\times\underline{s}}b_{kl}E_{kl}^{(2)}\right).
\end{equation}
Note that although $A$ and $B$ are matrices, the $E_{ij}^{(1)}\otimes E_{kl}^{(2)}=p_{ijkl}$ are  basis elements of $P$.  The matrix structure is not utilized here except in the indexing of coefficients.
 \end{srem}
 %SREM
 \begin{srem}[\bfseries Comments about basis notation, $E_i=\{e_{i1}, \ldots e_{in_i}\}$]
 \label{srem:nttstnbss}
 We have stated our standard assumptions for bases as follows:
``Let $V_1, \ldots, V_m$, $\dim(V_i) = n_i$, $i\in \underline{m}$, be vector spaces over $\bF$. Suppose, for all $i\in \underline{m}$, $V_i = \left<E_i\right>$,  
$E_i = \{e_{i1},\ldots, e_{in_i}\}$ is an ordered basis for $V_i$.''
Suppose, analogous to~\ref{rem:tnsprdmtr}, the 
$V_i = M_{p_i, q_i}$, $i=1, \ldots, m$, are vector spaces of matrices.
The standard ordered basis for $V_i$ is now 
\[
E_i=\{E_{11}^{(i)}, E_{12}^{(i)}, \ldots, E_{p_iq_i}^{(i)}\}.
\]
The basis element $e_{ij}$ has been replaced by $E_{st}^{(i)}$ where
$(s,t)$ is the element in position $j$ in the list 
$\underline{p_i}\times \underline{q_i}$ in lexicographic order.
We don't attempt to formalize this type of variation from the standard notation.
The general case will involve ordered bases where the order  is specified by a linear order on indices in some manner (usually some type of lex order).
 \end{srem}
 %SREM
 \begin{srem}[\bfseries Universal factorization property: example]
 \label{rem:unvfctprpexm}
This example continues the discussion of~\ref{rem:tnsprdmtr}
in order to illustrate lemma~\ref{lem:tnsrprdvctrsp}, the universal factorization property.
Let $V_1=\mathbf{M}_{p,q}$ ($p\times q$ matrices over $\bF$), and let $V_2=\mathbf{M}_{q,s}$. 
Let $A=(a_{ij})\in \mathbf{M}_{p,q}$ and $B=(b_{ij})\in \mathbf{M}_{q,s}$.
Here we have set $r=q$ in~\ref{rem:tnsprdmtr}.
Take the $\varphi$ of lemma~\ref{lem:tnsrprdvctrsp} to be the bilinear function
$\varphi(A,B) = AB$, the matrix product of $A$ and $B$.
Thus, $\varphi\in M(V_1, V_2: U)$ where $U =  \mathbf{M}_{p,s}.$
We want to construct the linear function 
\[
h_{\varphi}\in\bL(\mathbf{M}_{p,q}\otimes \mathbf{M}_{q,s}, \mathbf{M}_{p,s})
\] 
such that $h_{\varphi}(A\otimes B) = \varphi(A, B)$.
The set of matrices
$\{E^{(3)}_{cf}\Mid (c,f)\in \underline{p}\times\underline{s}\}$
is the standard basis for $U$  (notation as in~\ref{rem:tnsprdmtr}).
Following discussion~\ref{rem:vcspcfld}, we work with the component functions
 $\varphi^{(i,j)} \in M(V_1,V_2:\langle E^{(3)}_{ij}\rangle)$ for each $i,j$ fixed:
 \[
h_{\varphi^{(i,j)}}\in\bL(\mathbf{M}_{p,q}\otimes \mathbf{M}_{q,s}, 
\langle E^{(3)}_{ij}\rangle)\equiv 
\bL(\mathbf{M}_{p,q}\otimes \mathbf{M}_{q,s}, 
\bF).
\] 
Define $h_{\varphi^{(i,j)}}$ on the basis elements 
\[
\{E^{(1)}_{cd}\otimes E^{(2)}_{ef}\Mid (c,d)\in \underline{p}\times\underline{q},
(e,f)\in \underline{q}\times\underline{s}\}
\]
analogous to equation~\ref{eq:mltext2}: 
\begin{equation}
\label{eq:unvfctprpexm2}
h_{\varphi^{(i,j)}}(E^{(1)}_{cd}\otimes E^{(2)}_{ef}) = \varphi^{(i,j)}(E^{(1)}_{cd},E^{(2)}_{ef})=
E^{(1)}_{cd}E^{(2)}_{ef}(i,j)
\end{equation} 
where $E^{(1)}_{cd}E^{(2)}_{ef}$ denotes the matrix product of these two basis elements.
An elementary result from matrix theory states
$E^{(1)}_{cd}E^{(2)}_{ef}=\Ki(d=e)E^{(3)}_{cf}.$
Thus, 
\[
E^{(1)}_{cd}E^{(2)}_{ef}(i,j) =\Ki \left(d=e\right)E^{(3)}_{cf}(i,j)=
\Ki \left(d=e\right) \Ki \left((c,f)=(i,j)\right).
\]
Defining $h_{\varphi^{(i,j)}}$ on $A\otimes B$ by linear extension gives
\begin{equation}
\label{eq:unvfctprpexm3}
h_{\varphi^{(i,j)}}\left(A\otimes B\right) \coloneq 
\sum_{(c,d,e,f)} a_{cd}b_{ef}h_{\varphi^{(i,j)}}(E^{(1)}_{cd}\otimes E^{(2)}_{ef})
\end{equation}
Thus,
\begin{equation}
\label{eq:tnsprdtwomtr5}
h_{\varphi^{(i,j)}}\left(A\otimes B\right) = 
\sum_{(c,d,e,f)} a_{cd}b_{ef}\Ki \left(d=e\right) \Ki \left((c,f)=(i,j)\right).
\end{equation}
Thus,
\begin{equation}
\label{eq:tnsprdtwomtr6}
h_{\varphi^{(i,j)}}(A\otimes B) =\left( \sum_{d\in \underline{q}} a_{id}b_{dj}\right) = AB(i,j).
\end{equation}
Thus, 
\[
h_{\varphi} = \bigoplus_{(i,j)\in\Gamma(p,s)} h_{\varphi^{(i,j)}}
\] 
is in 
\[
\bL(\mathbf{M}_{p,q}\otimes \mathbf{M}_{q,s}, \mathbf{M}_{p,s})
\] 
and satisfies $h_{\varphi}(A\otimes B)=AB$.
\end{srem}

\begin{srem}[\bfseries Tensor products of matrices as matrices: a bad choice]
\label{rem:tnsprdmtr2}
Let $V_1=\mathbf{M}_{p_1,q_1}$, and let $V_2=\mathbf{M}_{p_2,q_2}$. We construct a tensor product $(P,\nu)$ for $V_1, V_2$. 
Let $n_1=p_1q_1$ and $n_2=p_2q_2$ and choose $P=\mathbf{M}_{n_1,n_2}$.
Define $\nu( E_{ij}^{(1)},E_{kl}^{(2)})\coloneq E_{ij}^{(1)}\otimes E_{kl}^{(2)} = p_{ijkl}$
where $p_{ijkl} = E^{(3)}_{pq}$ is the standard basis element of  
$\mathbf{M}_{n_1,n_2}$ with $p$ the position of $(i,j)$ in the lexicographic list
of $\underline{p_1}\times \underline{q_1}$ and $q$ defined similarly for 
$\underline{p_2}\times \underline{q_2}$.
By the construction of~\ref{rem:tnsprdalwexs} such a choice is possible. 
We have represented  $\mathbf{M}_{p_1,q_1}\otimes \mathbf{M}_{p_2,q_2}$ as matrices so that $A_1\otimes A_2$ is an $n_1\times n_2$ matrix where $n_1$ is the number of entries in $A_1$ and $n_2$ the number of entries in $A_2$.
The problem with this is that in the obvious extension to 
$V_1\otimes V_2\otimes V_3$ we would have
$\left(V_1\otimes V_2\right)\otimes V_3$ and $V_1\otimes \left(V_2\otimes V_3\right)$
isomorphic as tensor spaces (see~\ref{rem:asslwstnsprd}) but, in general, 
$(A_1\otimes A_2)\otimes A_3$ not equal to $A_1\otimes (A_2\otimes A_3)$
 as matrices.
(the former has $n_1n_2$ rows and the latter $n_1$ rows).
 \end{srem} 
 
%SREM
\begin{srem}[\bfseries Tensor products of matrices, a good choice: Kronecker product]
\label{rem:krnprd}
Let $V_1=\mathbf{M}_{p_1,q_1}$, and let $V_2=\mathbf{M}_{p_2,q_2}$. 
As in~\ref{rem:tnsprdmtr2}, we construct a tensor product $(P,\nu)$ for 
$V_1,  V_2$. 
Choose $P\coloneq\mathbf{M}_{p,q}$ where $p=p_1p_2$ and $q=q_1q_2$.
Let $\mu\in \underline{p_1}\times \underline{p_2}$ and 
$\kappa\in \underline{q_1}\times \underline{q_2}$ (both with lexicographic order).
We use these two ordered sets to index the rows and columns of the 
matrices in  $P$.  
The matrices 
\begin{equation}
\label{eq:krnprd0}
\{E_{\mu,\kappa} \Mid \mu\in \underline{p_1}\times \underline{p_2}\,,\,  \kappa\in \underline{q_1}\times \underline{q_2}\,\}
\end{equation}
where $E_{\mu,\kappa}(\alpha, \beta)=\mathcal{X}((\alpha,\beta)=(\mu,\kappa))$
are the standard basis elements of $\mathbf{M}_{p,q}$. 
Let 
\[
\{E^{(1)}_{\mu(1),\kappa(1)}\Mid (\mu(1),\kappa(1))\in\underline{p_1}\times \underline{q_1}\}
\] 
and 
\[\{E^{(2)}_{\mu(2),\kappa(2)}\Mid (\mu(2),\kappa(2))\in\underline{p_2}\times \underline{q_2}\}, 
\]
each ordered lexicographically, denote the standard ordered bases for $\mathbf{M}_{p_1,q_1}$ and  
$V_2=\mathbf{M}_{p_2,q_2}$.

Define (see~\ref{lem:cnnbss}) the function $\nu$ of $(P,\nu)$ by
\begin{equation}
\label{eq:krnprd1}
\nu(E^{(1)}_{\mu(1),\kappa(1)}, \,E^{(2)}_{\mu(2),\kappa(2)})\coloneq 
E^{(1)}_{\mu(1),\kappa(1)} \otimes E^{(2)}_{\mu(2),\kappa(2)} = E_{\mu,\kappa}.
\end{equation}
Let $A_i=(a_i(s,t))\in \mathbf{M}_{p_i,q_i}$, $i=1,2$. Applying multilinear extension~\ref{eq:mltlnrfncb}, we get $\nu(A_1,A_2)=A_1\otimes A_2=$
\begin{equation}
\label{eq:krnprd2}
\nu\left(\sum_{\mu(1),\kappa(1)}a_1(\mu(1),\kappa(1))E^{(1)}_{\mu(1),\kappa(1)},
\sum_{\mu(2),\kappa(2)}a_2(\mu(2),\kappa(2))E^{(2)}_{\mu(2),\kappa(2)}\right)=
\end{equation}
\[
\sum_{\mu,\kappa} \prod_{i=1}^2 a_i(\mu(i),\kappa(i)) E_{\mu,\kappa}.
\]
Equation~\ref{eq:krnprd1} guarantees  that $(P,\nu)$ is a tensor 
product of $V_1$ and $V_2$ (\ref{rem:eqvdfnrtnsprd}). 
The notation of~\ref{eq:krnprd2} is chosen to extend from $m=2$ to the general case. 
The tensor product, $(P,\nu)$, $P=\otimes_{i=1}^2\mathbf{M}_{p_i,q_i}$ and $\nu$ defined by equation~\ref{eq:krnprd1} and equation~\ref{eq:krnprd2} is called the {\em Kronecker tensor product} of 
the vector spaces $\mathbf{M}_{p_1,q_1}$ and $\mathbf{M}_{p_2,q_2}.$
From equation~\ref{eq:krnprd2}, we see that the homogeneous tensors,
$\nu(A_1,A_2)\coloneq A_1\otimes A_2$ can be regarded as matrices in 
$\mathbf{M}_{p,q}$ and that the $\mu, \kappa$ entry of $A_1\otimes A_2$ is
\begin{equation}
\label{eq:krnprd3}
\left(A_1\otimes A_2\right)(\mu, \kappa) = 
\prod_{i=1}^2 A_i(\mu(i),\kappa(i)).
\end{equation}
These homogeneous elements are are sometimes defined without reference  to the general theory of tensor spaces  and are called  ``Kronecker products of matrices.''
\end{srem}

\begin{srem}[\bfseries Kronecker product: example of homogeneous tensors]
\label{srem:krnprdexm} 
For
\[A_1=(a_1(i,j)),\; A_2=(a_2(i,j))\in \mathbf{M}_{2,2}\] 
the Kronecker product
$A_1\otimes A_2 = $
\begin{equation*} 
\bordermatrix{~ &1\,1 & 1\,2 &2\,1&2\,2 \cr
%L1
\underline{1}\,\underline{1}& a_1(\underline{1},1)a_2(\underline{1},1) & a_1(\underline{1},1)a_2(\underline{1},2) & a_1(\underline{1},2)a_2(\underline{1},1) & a_1(\underline{1},2)a_2(\underline{1},2) \cr
%L2
\underline{1}\,\underline{2}& a_1(\underline{1},1)a_2(\underline{2},1) & a_1(\underline{1},1)a_2(\underline{2},2) & a_1(\underline{1},2)a_2(\underline{2},1) & a_1(\underline{1},2)a_2(\underline{2},2) \cr
%L3
\underline{2}\,\underline{1}& a_1(\underline{2},1)a_2(\underline{1},1) & a_1(\underline{2},1)a_2(\underline{1},2) & a_1(\underline{2},2)a_2(\underline{1},1) & a_1(\underline{2},2)a_2(\underline{1},2) \cr
%L4 
\underline{2}\,\underline{2}& a_1(\underline{2},1)a_2(\underline{2},1) & a_1(\underline{2},1)a_2(\underline{2},2) & a_1(\underline{2},2)a_2(\underline{2},1) & a_1(\underline{2},2)a_2(\underline{2},2) \cr
}
\end{equation*}

Note that  $A_1\otimes A_2$ can be constructed by starting with a  $2\times 2$ matrix with entries in  $\mathbf{M}_{2,2}$,
\[
\begin{pmatrix}
a_1(1,1)A_2\; & \;a_1(1,2)A_2\\
a_1(2,1)A_2\; & \;a_1(2,2)A_2
\end{pmatrix},
\]
and doing the indicated multiplications of entries from $A_1$ with $A_2$ to construct the blocks of $A_1\otimes A_2$.  
This pleasing structure is a consequence of using lexicographic order as done in~\ref{rem:krnprd}.
Other orders would work just as well, but the result might be a mess to the human eye.
\end{srem}
%SREM
\begin{srem}[\bfseries General Kronecker products of vector spaces of matrices]
\label{srem:krnprdgnrcas}
We follow the discussion~\ref{rem:krnprd}, developing the notation for the general case.
Define a tensor product $(P,\nu)$ for $V_1,\ldots, V_m$ where
 $V_1=\mathbf{M}_{p_1,q_1},\ldots,V_m=\mathbf{M}_{p_m,q_m}$. 
 Choose $P\coloneq\mathbf{M}_{p,q}$ where $p=\prod_{i=1}^mp_i$ and 
$q=\prod_{i=1}^mq_i$.
Let $\mu\in \underline{p_1}\times\cdots\times \underline{p_m}$ and 
$\kappa\in \underline{q_1}\times\cdots\times \underline{q_m}$ (both with lexicographic order).
We use these two ordered sets to index, respectively, the rows and columns of the matrices in  $P$.  
We choose the matrices 
\begin{equation}
\label{eq:krnprd1g}
\{E_{\mu,\kappa} \Mid \mu\in \underline{p_1}\times\cdots\times \underline{p_m}\,,\, 
\kappa\in \underline{q_1}\times\cdots\times \underline{q_m}\,\}
\end{equation}
where $E_{\mu,\kappa}(\alpha, \beta)=\mathcal{X}((\alpha,\beta)=(\mu,\kappa))$
to be the basis elements of $\mathbf{M}_{p,q}$. 
Let
\[
\{E^{(i)}_{\mu(i),\kappa(i)}\Mid (\mu(i),\kappa(i))\in\underline{p_i}\times \underline{q_i}\},
\;i=1,\ldots, m,
\] 
be the standard bases for the $V_i=\mathbf{M}_{p_i,q_i}$.
We define $\nu \in M(V_1, \ldots, V_m: P)$ by
\begin{equation}
\label{eq:krnprd1gnr}
\nu(E^{(1)}_{\mu(1),\kappa(1)},\ldots, \,E^{(m)}_{\mu(m),\kappa(m)})\coloneq 
E^{(1)}_{\mu(1),\kappa(1)} \otimes\cdots\otimes E^{(m)}_{\mu(2),\kappa(m)} 
= E_{\mu,\kappa}.
\end{equation}
Let $A_i=(a_i(s,t))\in \mathbf{M}_{p_i,q_i}$, $i=1,\ldots, m$. 
Applying multilinear extension~\ref{eq:mltlnrfncb}, we get
 $\nu(A_1,\ldots, A_m)=A_1\otimes\cdots\otimes A_m=$
\begin{equation}
\label{eq:krnprd2gnr}
\nu\left(\sum_{\mu(1),\kappa(1)}a_1(\mu(1),\kappa(1))E^{(1)}_{\mu(1),\kappa(1)},
\ldots,
\sum_{\mu(m),\kappa(m)}a_m(\mu(m),\kappa(m))E^{(m)}_{\mu(m),\kappa(m}\right)=
\end{equation}
\[
\sum_{\mu,\kappa} \prod_{i=1}^m a_i(\mu(i),\kappa(i)) E_{\mu,\kappa}.
\]
\end{srem}
%DEF
\begin{definition}[\bfseries Kronecker product general definition]
\label{def:krnprdgnrdfn}
Consider the vector spaces of matrices $\mathbf{M}_{p_1,q_1},\ldots, \mathbf{M}_{p_m,q_m}$ and  $\mathbf{M}_{p,q}$ 
where $p=\prod_{i=1}^m p_i$, $q=\prod_{i=1}^m q_i$, $i=1,\ldots, m$. 
Order the rows of the matrices in $\mathbf{M}_{p,q}$ with the 
$
\{\mu\Mid \mu\in \underline{p_1}\times\cdots\times \underline{p_m}\}
$ 
and the columns with
$
\{\kappa\Mid\kappa\in \underline{q_1}\times\cdots\times \underline{q_m}\},
$
each set ordered lexicograhically.
The pair $(\mathbf{M}_{p,q}, \nu )$ is a tensor product of these vector spaces where
$p=\prod_{i=1}^m p_i$, $q=\prod_{i=1}^m q_i$, $i=1,\ldots, m$, and 
$\nu\in M(\mathbf{M}_{p_1,q_1}, \ldots, \mathbf{M}_{p_m,q_m}: \mathbf{M}_{p,q})$
is defined by $\nu(A_1, \ldots, A_m)\coloneq A_1\otimes\cdots \otimes A_m$ where
\[
(A_1\otimes\cdots \otimes A_m)(\mu, \kappa) = \prod_{i=1}^m  A_i(\mu(i),\kappa(i)). 
\]
\end{definition}
\begin{srem}[\bfseries Is the Kronecker tensor product a special case?]
\label{srem:krnprdgnrdsc}
In ~\ref{rem:eqvdfnrtnsprd} we give three equivalent definitions of a tensor product of vector spaces.
The setup for these definitions is as follows:

``Let $V_1, \ldots, V_m$, $\dim(V_i) = n_i$, be vector spaces over $\bF$ . 
Let $P$ be a vector space over $\bF$ and $\nu\in M(V_1, \ldots, V_m:P).$
Suppose, for all $i=1,\ldots m$,  
$E_i = \{e_{i1},\ldots, e_{in_i}\}$ is a basis for $V_i$.''

The Kronecker product is usually described as a ``special case'' of a tensor product.
It is actually equivalent to the definition of a tensor product, differing only by the indexing 
of the bases of the $V_i$ and $P$.  If we factor the dimensions $n_i=p_iq_i$ (possible in many ways; we allow $p_i=1$ or $q_i=1$), then we can replace $E_i = \{e_{i1},\ldots, e_{in_i}\}$ by
$E^{(i)}=\{E^{(i)}_{11}, \ldots, E^{(i)}_{p_iq_i}\}$ ordered, for example, lexicographically
on indices.   The vector space $P$ has dimension $\prod_{i=1}^m n_i$ which is the same
as $pq$ where $p=\prod_{i=1}^m p_i$ and $q=\prod_{i=1}^m q_i$. 
The most direct correspondence is  to replace 
$E_i = \{e_{i1},\ldots, e_{in_i}\}$ 
by
$E^{(i)}=\{E^{(i)}_{11}, \ldots, E^{(i)}_{1,n_i}\}$.
The matrix model for Kronecker products is important for applications to matrix theory that are model specific.
\end{srem}
\section{Inner products on tensor spaces}
\label{sec:innerproducts}
%SREM
\begin{srem}[\bfseries Inner products and tensor spaces]
\label{srem:innprdtnsrspcx}
Assume $V_1, \ldots, V_m$, $\dim(V_i) = n_i$,   are vector spaces over the complex 
numbers, $\bC$. Let $(W, \omega)$ be a tensor product of these vector spaces.
As usual, $W=\otimes_{i=1}^mV_i$, and 
$\omega(x_1, \ldots, x_m)=x_1{\otimes}\cdots\otimes x_m$ are the homogeneous tensors.
Recall that  $\varphi\in M(W,W:\bC)$ is {\em conjugate bilinear} if for all 
$a,a',b, b'\in W$
 and $c,d\in \bC$ 
\begin{equation}
\label{eq:cnjblncnd}
\begin{matrix}
\varphi(ca+da',b)=c\varphi(a,b) + d\varphi(a',b) \\
\varphi(b,ca+da')=\overline{c}\varphi(b,a) + \overline{d}\varphi(b,a').
\end{matrix}
\end{equation}
A conjugate bilinear $\varphi\in M(W,W:\bC)$ is an {\em inner product} on $W$ 
if for all $a, b\in W$ 
\begin{equation}
\label{eq:innprdunt}
\begin{matrix}
\varphi(a,b)=\overline{\varphi(b,a)}\\
\varphi(a,a)\geq 0\;\mathrm{with}\; \varphi(a,a) = 0 \;\mathrm{iff}\; a=0.
\end{matrix}
\end{equation}
In that case, the pair $(W, \varphi)$ is  a {\em unitary space}.
The first condition of equation~\ref{eq:innprdunt} is called {\em conjugate symmetric}, the second is called {\em positive definite}.

If $\{f_{\gamma}\Mid \gamma\in \Gamma(n_1, \ldots, n_m)\}$ is any 
basis of $W$ then an inner product on $W$ can be defined by specifying
$\varphi(f_{\alpha}, f_{\beta}) = \delta_{\alpha, \beta}$ and 
defining $\varphi\in M(W,W:\bC)$ by extending these values by (conjugate) bilinear extension.
In this case, $\{f_{\gamma}\Mid \gamma\in \Gamma(n_1, \ldots, n_m)\}$ is called an 
{\em orthonormal basis} for the inner product $\varphi$.
If $\varphi$ is any inner product on $W$, there exists an orthonormal basis that defines 
$\varphi$ in the manner just described (e.g., by using the Gram-Schmidt orthonormalization process).

If $\{f_{\gamma}\Mid \gamma\in \Gamma(n_1, \ldots, n_m)\}$ is an {\em orthonormal basis} for the inner product $\varphi$ and 
$a=\sum_{\alpha\in \Gamma}c_{\alpha}f_{\alpha}$ and 
$b=\sum_{\beta\in \Gamma}d_{\beta}f_{\beta}$ then
\begin{equation}
\label{eq:gnrinnprdsm0}
\varphi(a,b) = 
\sum_{\alpha\in \Gamma}\sum_{\beta\in \Gamma}c_{\alpha}\overline{d_{\beta}}
 \varphi(f_{\alpha},f_{\beta}) = 
 \sum_{\alpha\in \Gamma}c_{\alpha}\overline{d_{\alpha}}.
\end{equation}
We want to relate these ideas more closely to the tensor product, $(W, \omega)$.
Suppose, for $i=1,\ldots m$,  
$E_i = \{e_{i1},\ldots, e_{in_i}\}$ is an orthonormal basis for $V_i$ with inner 
product $\varphi_i \in M(V_i,V_i: \bC)$.
Assume $x_i = \sum_{k=1}^{n_i} c_{ik} e_{ik}$ and assume 
$y_i = \sum_{k=1}^{{n}_i} d_{ik} {e}_{ik}.$ 
We use the notation 
\begin{equation}
\label{eq:shrtnttxotimes0}
x^{\otimes} \coloneq x_1\otimes\cdots\otimes x_m\;\mathrm{and}\;
e_{\gamma}^{\otimes} \coloneq e_{1\gamma(1)}\otimes\cdots \otimes e_{m\gamma(m)}.
\end{equation}
\end{srem}
%LEM
\begin{lemma}[\bfseries Inner product on $W$ as a product]
\label{lem:innprdasprd}
We refer to~\ref{srem:innprdtnsrspcx} for notation.
Suppose, for $i=1,\ldots m$,  
$E_i = \{e_{i1},\ldots, e_{in_i}\}$ is an orthonormal basis for $V_i$ with inner 
product $\varphi_i \in M(V_i,V_i: \bC)$.
There exists a unique inner product 
$\varphi\in M(W,W:\bC)$, $W=\otimes_{i=1}^m V_i$, such that
\begin{equation}
\label{eq:cnjblnfncprdfrm0}
\varphi(x^{\otimes}, y^{\otimes})=\varphi (x_1\otimes\cdots\otimes x_m, y_1\otimes\cdots\otimes y_m) =
\prod_{i=1}^m \varphi_i(x_i,y_i).
\end{equation}
This $\varphi$ is defined  by  
$\varphi(e_{\alpha}^{\otimes},e_{\beta}^{\otimes}) \coloneq
\prod_{i=1}^m \varphi_{i}(e_{i\alpha(i)}, e_{i\beta(i)})$ for all 
$\alpha, \beta\in \Gamma(n_1,\ldots, n_m)$.  The basis
$\{e_{\alpha}^{\otimes}\Mid \alpha \in \Gamma(n_1,\ldots, n_m)\}$ is an orthonormal basis for $\varphi$.
\begin{proof}
Define $\varphi(e_{\alpha}^{\otimes},e_{\beta}^{\otimes}) \coloneq
\prod_{i=1}^m \varphi_{i}(e_{i\alpha(i)}, e_{i\beta(i)}),\,$ 
$\alpha, \beta\in \Gamma(n_1,\ldots, n_m)$ and define 
$\varphi\in M(W,W:\bC)$ by conjugate bilinear extension.
Note that $\prod_{i=1}^m \varphi_{i}(e_{i\alpha(i)}, e_{i\beta(i)})=\delta_{\alpha\beta}$
so that the $\varphi$ so defined has   
$\{e_{\alpha}^{\otimes}\Mid \alpha \in \Gamma(n_1,\ldots, n_m)\}$
as an orthonormal basis.
We have
\[
\prod_{i=1}^m \varphi_i(x_i,y_i) = 
\prod_{i=1}^m \varphi_i\left( \sum_{j=1}^{n_i} c_{ij} e_{ik},
\sum_{k=1}^{n_i} d_{ik} e_{ik}\right)=
\]
\[
\prod_{i=1}^m \left(\sum_{j, k} c_{ij}\overline{d_{ik}} \varphi_i(e_{ij}, e_{ik})\right) =
\prod_{i=1}^m \left(\sum_{t=1}^{n_i} c_{it}\overline{d_{it}}\right) =
\]
\[
\sum_{\alpha\in \Gamma} \prod_{i=1}^{m} c_{i\alpha(i)}\overline{d_{i\alpha(i)}}
=\sum_{\alpha\in \Gamma} \prod_{i=1}^{m} c_{i\alpha(i)}\prod_{i=1}^{m} \overline{d_{i\alpha(i)}}=
\varphi (x^{\otimes}, y^{\otimes}).
\]
The last equality follows from equation~\ref{eq:gnrinnprdsm0} with
$a=x^{\otimes}$ and $b= y^{\otimes}$ noting that 
$c_{\alpha}=\prod_{i=1}^{m} c_{i\alpha(i)}$ and 
$\overline{d_{\alpha}}=\prod_{i=1}^{m} \overline{d_{i\alpha(i)}}$ for these choices.
\end{proof}
\end{lemma}
%SREM
We now extend lemma~\ref{lem:innprdasprd} to the case of conjugate bilinear functions that need not be inner products.
\begin{srem}[\bfseries Conjugate bilinear functions on tensor spaces, general remarks]
\label{srem:ecluntinnprd}
We specify some notational conventions to be used in lemma~\ref{lem:indinprd2}. Assume $V_1, \ldots, V_m$, $\dim(V_i) = n_i$, and 
$\hat{V}_1, \ldots, \hat{V}_m$, $\dim(\hat{V}_i) = \hat{n}_i$, are vector spaces over the complex numbers, $\bC$. 
Let $W=\otimes_{i=1}^mV_i$ and   $\hat{W}=\otimes_{i=1}^m\hat{V}_i.$
Suppose, for $i=1,\ldots m$,  
$E_i = \{e_{i1},\ldots, e_{in_i}\}$ is a basis for $V_i$ and
$\hat{E}_i = \{\hat{e}_{i1},\ldots, \hat{e}_{in_i}\}$ is a basis for $\hat{V}_i.$
Assume $x_i = \sum_{k=1}^{n_i} c_{ik} e_{ik}$ and  $y_i = \sum_{k=1}^{\hat{n}_i} d_{ik} \hat{e}_{ik}.$ 
We use the notation 
\begin{equation}
\label{eq:shrtnttxotimes}
x^{\otimes} \coloneq x_1\otimes\cdots\otimes x_m\;\mathrm{and}\;
e_{\gamma}^{\otimes} \coloneq e_{1\gamma(1)}\otimes\cdots \otimes e_{m\gamma(m)}.
\end{equation}
Thus, $\hat{e}_{\gamma}^{\otimes}=\hat{e}_{1\gamma(1)}\otimes\cdots \otimes \hat{e}_{m\gamma(m)}.$
Let $\Gamma = \Gamma(n_1,\ldots, n_m)$ and  
$\hat{\Gamma} = \Gamma(\hat{n}_1,\ldots, \hat{n}_m)$
(notation~\ref{srem:tnsprdsts}).
Let $\{e_{\gamma}^{\otimes}\Mid \gamma\in \Gamma\}$ be the basis of $W$ induced by the bases $E_i$, and define $\{\hat{e}_{\gamma}^{\otimes}\Mid \gamma\in \hat{\Gamma\}}$ 
similarly for $\hat{E}_i$ and $\hat{W}.$
Assume that $\varphi\in M(W,\hat{W}:\bC)$ is conjugate bilinear: for all $a,a',b, b'\in W$
or $\hat{W}$ (as appropriate) and $c,d\in \bC$
\begin{equation}
\label{eq:cnjblncnd2}
\begin{matrix}
\varphi(ca+da',b)=c\varphi(a,b) + d\varphi(a',b) \\
\varphi(b,ca+da')=\overline{c}\varphi(b,a) + \overline{d}\varphi(b,a').
\end{matrix}
\end{equation}

If $\{f_{\gamma}\Mid \gamma\in \Gamma\}$, $\{\hat{f}_{\gamma}\Mid \gamma\in \hat{\Gamma\}}$ 
are bases for $W$, $\hat{W}$, respectivly and 
$a=\sum_{\alpha\in \Gamma}c_{\alpha}f_{\alpha}$ and 
$b=\sum_{\beta\in \hat{\Gamma}}d_{\beta}\hat{f}_{\beta}$ then
\begin{equation}
\label{eq:gnrinnprdsm}
\varphi(a,b) =  
\sum_{\alpha\in \Gamma}\sum_{\beta\in \hat{\Gamma}}c_{\alpha}\overline{d_{\beta}}
 \varphi(f_{\alpha},\hat{f}_{\beta}).
\end{equation}
\end{srem}

\begin{lemma}[\bfseries Conjugate bilinear functions as products]
\label{lem:indinprd2}
We use the terminology of~\ref{srem:ecluntinnprd}.
Assume $\varphi_i\in M(V_i,\hat{V}_i:\bC)$, $i\in\underline{m}$, is conjugate bilinear.  
There exists a unique conjugate bilinear function 
$\varphi\in M(W,\hat{W}:\bC)$ such that
\begin{equation}
\label{eq:cnjblnfncprdfrm1}
\varphi (x_1\otimes\cdots\otimes x_m, y_1\otimes\cdots\otimes y_m) =
\prod_{i=1}^m \varphi_i(x_i,y_i).
\end{equation}
Define $\varphi$ by  
$\varphi(e_{\alpha}^{\otimes},\hat{e}_{\beta}^{\otimes}) \coloneq
\prod_{i=1}^m \varphi_{i}(e_{i\alpha(i)}, \hat{e}_{i\beta(i)})$ for all 
$\alpha\in \Gamma(n_1,\ldots, n_m)$ and 
$\beta\in \Gamma(\hat{n}_1,\ldots, \hat{n}_m)$.
\begin{proof}
In equation~\ref{eq:gnrinnprdsm}, take 
$f_{\gamma} := e_{\gamma}^{\otimes}$, $\gamma\in \Gamma$ and
$\hat{f}_{\gamma} := \hat{e}_{\gamma}^{\otimes}$, $\gamma\in \hat{\Gamma}$, and
let $a=x^{\otimes}$, $b=y^{\otimes}$ (~\ref{srem:ecluntinnprd}).
We have  $x^{\otimes}= \sum_{\alpha\in \Gamma}(\prod_{i=1}^m c_{i\alpha(i)})e_{\alpha}^{\otimes}$
and $y^{\otimes}= \sum_{\beta\in\hat{\Gamma}}(\prod_{i=1}^m d_{i\beta(i)})\hat{e}_{\beta}^{\otimes}.$
Equation~\ref{eq:gnrinnprdsm} becomes
\begin{equation}
\label{eq:gnrinnprdsmvar}
\varphi(x^{\otimes},y^{\otimes}) = \sum_{\alpha\in \Gamma}\sum_{\beta\in {\hat{\Gamma}}}
\prod_{i=1}^m c_{i\alpha(i)}\prod_{i=1}^m \overline{d}_{i\beta(i)}
\varphi(e_{\alpha}^{\otimes},\hat{e}_{\beta}^{\otimes}).
\end{equation}
Define $\varphi(e_{\alpha}^{\otimes},\hat{e}_{\beta}^{\otimes}) \coloneq 
\prod_{i=1}^m \varphi_{i}(e_{i\alpha(i)}, \hat{e}_{i\beta(i)}).$ This specifies $\varphi$ uniquely as a 
conjugate bilinear function on $W$ (equation~\ref{eq:gnrinnprdsm}).
We need to show that~\ref{eq:cnjblnfncprdfrm1} holds.
Write~\ref{eq:gnrinnprdsmvar} as 
\begin{equation}
\label{eq:gnrinnprdsmvar2}
\varphi(x^{\otimes},y^{\otimes}) = 
\sum_{\alpha\in \Gamma}\prod_{i=1}^m c_{i\alpha(i)}
\left(\sum_{\beta\in \hat{\Gamma}}\prod_{i=1}^m 
\overline{d}_{i\beta(i)}\varphi_{i}(e_{i\alpha(i)}, \hat{e}_{i\beta(i)})\right).
\end{equation}
Next, we interchange product and sum (\ref{eq:changesumprod}) on the expression inside the parentheses
of~\ref{eq:gnrinnprdsmvar2}.  With 
$a_{ij} =\overline{d}_{ij}\varphi_{i}(e_{i\alpha(i)}, \hat{e}_{ij})$ we obtain
\begin{equation}
\label{eq:gnrinnprdsmvar3}
\varphi(x^{\otimes},y^{\otimes}) = 
\sum_{\alpha\in \Gamma}\prod_{i=1}^m c_{i\alpha(i)}
\prod_{i=1}^m\left(\sum_{j=1}^{n_i} 
\overline{d}_{ij}\varphi_{i}(e_{i\alpha(i)}, \hat{e}_{ij})\right).
\end{equation}
Using that the $\varphi_i$ are conjugate bilinear and recalling that 
$y_i=\sum_{j=1}^{n_i} d_{ij}\hat{e}_{ij}$, we obtain
\begin{equation}
\label{eq:gnrinnprdsmvar4}
\varphi(x^{\otimes},y^{\otimes}) = 
\sum_{\alpha\in \Gamma}\prod_{i=1}^m c_{i\alpha(i)}\varphi_{i}(e_{i\alpha(i)}, y_i).
\end{equation}
Repeating this sum-interchange process with $a_{ij}=c_{ij}\varphi_{i}(e_{ij}, y_i)$ gives
\begin{equation}
\label{eq:gnrinnprdsmvar5}
\varphi(x^{\otimes},y^{\otimes}) = 
\prod_{i=1}^m \varphi_{i}(x_i, y_i).
\end{equation}
By definition, $\varphi_i(x_i,y_i)$ are independent of the $E_i$ and $\hat{E}_i$, $i\in\underline{m}.$
\end{proof}
\end{lemma}
%SREM
%\begin{srem}[\bfseries Example of conjugate  bilinear functions as products]
%\label{srem:explcnjgblnrfuncprds}
%We follow the discussion~\ref{srem:ecluntinnprd} and apply 
%lemma~\ref{lem:indinprd2}.
%Let $V_1=\bC^2$, $V_2=\bC^3$, $\hat{V}_1=\bC^4$, and $\hat{V}_2=\bC^4$.
%Define $\varphi_1\in M(V_1, \hat{V}_1)$ by 
%\[
%\varphi_1((c_{11}, c_{12}), (d_{11}, d_{12}, d_{13}, d_{14})) = 
%\left(\sum_{i=1}^2 c_{1i}\right)\left(\sum_{i=1}^4 \overline{d_{1i}}\right).
%\]
%\end{srem}

\section{Direct sums and tensor products}
\label{sec:directsums}
%SREM PARTIONING
\begin{srem}[\bfseries Partitioning $\Gamma(n_1, \ldots, n_m)$]
\label{srem:prtngmma}
Let $\bD_i = \{D_{i1}, D_{i2}, \ldots D_{ir_i}\}$ be a partition of $\underline{n_i}$ for $i=1, 2, \ldots, m$ with blocks ordered as indicated.  
For every $\gamma\in\Gamma(n_1, \ldots, n_m)$ there exists a unique
$\alpha\in\underline{r_1}\times\cdots\times\underline{r_m}$ such that
$\gamma(i)\in D_{i\alpha(i)}$ for $i=1, 2, \ldots, m.$

Let $\bD_{\alpha}\coloneq \times_{i=1}^m D_{i\alpha(i)} = 
\{\gamma\Mid \gamma(i)\in D_{i\alpha(i)},\, i\in\underline{m}\}.
$ 
The set 
$\bD_{\Gamma}=\{\bD_{\alpha}\Mid \alpha\in\underline{r_1}\times\cdots\underline{r_m}\}$,  
is a partition of  $\Gamma(n_1, \ldots, n_m)$ with blocks 
$\bD_{\alpha}$.
At this point, we haven't specified an order on the blocks of $\bD_{\Gamma}$. 
A natural choice would be to use lexicographic order on the domain,
$\Gamma(r_1, \ldots, r_m)$, of the map $\alpha\mapsto \bD_{\alpha}$.

We introduce the following definition.
\end{srem}
%DEF
\begin{definition}[\bfseries Partition of $\Gamma$ induced by partitions of the $\underline{n_i}$]
\label{def:prtindnsb}
Let $\bD_i = \{D_{i1}, D_{i2}, \ldots D_{ir_i}\}$ be a partition of $\underline{n_i}$ for $i=1, \ldots, m$ with blocks ordered as indicated. 
The partition
\[ 
\bD_{\Gamma}=\{\bD_{\alpha}\Mid\bD_{\alpha}=\times_{i=1}^m D_{i\alpha(i)},\; \alpha\in\underline{r_1}\times\cdots\underline{r_m}\}
\]
is the {\em partition of}\hspace{3 pt} $\Gamma(n_1, \ldots, n_m)$ {\em induced by the 
partitions $\bD_i$, $i=1, \ldots, m.$}
\end{definition}

%SREM
\begin{srem}[\bfseries Examples of partitions of type $\bD_{\Gamma}$]
\label{srem:dscprtgmma}
Let $V_1$ and $V_2$ be vectors spaces with ordered bases
$E_1=\{e_{11}, e_{12}, e_{13}\}$ and $E_2=\{e_{21}, e_{22}, e_{23}, e_{24}\}$ respectively.
For convenience, we list the elements of $\Gamma(3,4)$ in lexicographic order ($m=2$ in this case): 
\[
(1,1), (1,2), (1,3), (1,4), (2,1), (2,2), (2,3), (2,4), (3,1), (3,2), (3,3), (3,4).\]
Take the ordered partitions $D_i$ to be as follows: $D_1 = \{D_{11}, D_{12}\}$ with  
$D_{11}=\{1,3\}$ and $D_{12}=\{2\}$;  $D_2 = \{D_{21}, D_{22}\}$ with 
$D_{21}=\{2,4\}$ and $D_{22}=\{1,3\}$.
\[ 
\bD_{\Gamma}=\{\bD_{\alpha}\Mid\bD_{\alpha}=\times_{i=1}^2 D_{i\alpha(i)},\; \alpha\in\underline{2}\times\underline{2}\}.
\]
We construct $\bD_{\alpha}$:

With $\alpha=(1,1)$ we get
\[
\bD_{(11)} = D_{11}\times D_{21} = \{1,3\}\times \{2,4\}= \{(1,2), (1,4), (3,2), (3,4)\}.
\]
With $\alpha = (1,2)$ we get
\[
\bD_{(12)} = D_{11}\times D_{22} = \{1,3\}\times \{1,3\}=\{(1,1), (1,3), (3,1), (3,3)\}.
\]
With $\alpha = (2,1)$ we get
\[
\bD_{(21)} = D_{12}\times D_{21} = \{2\}\times \{2,4\}=\{(2,2), (2,4)\}.
\]
With $\alpha = (2,2)$ we get
\[
\bD_{(22)} = D_{12}\times D_{22} = \{2\}\times \{1,3\}=\{(2,1), (2,3)\}.
\]
 Thus, 
 \[
 \bD_{\Gamma}=\{\{(1,2), (1,4), (3,2), (3,4)\}, \{(1,1), (1,3), (3,1), (3,3)\},
 \]
 \[
 \{(2,2), (2,4)\}, \{(2,1), (2,3)\}.
 \]
 The blocks of $\bD_{\Gamma}$ can be ordered by   lexicographically ordering the domain, $\Gamma({2}, {2})$,  of $\alpha\mapsto \bD_{\alpha}$ to get
$\bD_{(11)}, \bD_{(12)}, \bD_{(21)}, \bD_{(22)}.$
 If we let $W_{ij} = \langle e_{it}\Mid t\in D_{ij}\rangle$ then $V_i = \oplus_{j=1}^2 W_{ij}$.
 For example, $W_{21} = \langle e_{22}, e_{24}\rangle$ and 
 $W_{22} = \langle e_{21}, e_{23} \rangle$, and $V_2 = W_{21}\oplus W_{22}$. 
\end{srem}
\begin{remark}[\bfseries Direct sums and bases]
\label{rem:drcsmsbss}
The following is a slight generalization of a standard theorem from linear algebra:

Let $V_i$, $\dim(V_i)=n_i$, $i=1,\ldots, n_i$.
Let $\bD_i = \{D_{i1}, D_{i2}, \ldots D_{ir_i}\}$ be an ordered partition of 
$\underline{n_i}$ for $i=1, \ldots, m$.
Then $V_i = \oplus_{j=1}^{r_i} W_{ij}$, $\dim(W_{ij})=|D_{ij}|$, if and only if
there exists bases $e_{i1}, \ldots, e_{in_i}$ such that 
$W_{ij}=\langle e_{it}\Mid t\in D_{ij}\rangle$.
\end{remark} 
%DEF
%SREM
\begin{srem}[\bfseries Tensor products of direct sums]
\label{srem:tnsprddrsms}
Let $(P, \nu)$ be a tensor product (\ref{def:tnsrprdfld}) of $V_i$, 
$\dim(V_i)=n_i$ Let $E_i=\{e_{i1}, \ldots, e_{in_i}\}$ be ordered bases for the $V_i$, $i=1, \ldots, m$.
Let 
$\{p_{\gamma}\Mid \gamma\in \Gamma(n_1,\ldots, n_m), p_{\gamma} =\nu(e_{1\gamma(1)}, \ldots, e_{m\gamma(m)})\}$ be the associated basis.
 As in definition~\ref{def:prtindnsb}, examples~\ref{srem:dscprtgmma} and
 remark~\ref{rem:drcsmsbss},
let $\bD_i = \{D_{i1}, D_{i2}, \ldots D_{ir_i}\}$ be an ordered partition of $\underline{n_i}$ for 
$i=1, \ldots, m$. 
If $W_{ij} = \langle e_{it}\Mid t\in D_{ij} \rangle$ then $V_i = \oplus_{j=1}^{r_i} W_{ij}$,
$i=1\ldots m$.
From definition~\ref{def:prtindnsb} we have 
\[ 
\bD_{\Gamma}=\{\bD_{\alpha}\Mid\bD_{\alpha}=\times_{i=1}^m D_{i\alpha(i)},\; \alpha\in\underline{r_1}\times\cdots\times\underline{r_m}\}.
\]
Consider $\times_{i=1}^m W_{i\alpha(i)}$.
Let  $P_{\alpha} = \langle p_{\gamma}\Mid \gamma\in \bD_{\alpha}\rangle$.  
From the fact that $\bD_{\Gamma}$ is a partition of $\Gamma$, we have
$P=\oplus_{\alpha} P_{\alpha}$.
We claim that $(P_{\alpha}, \nu_{\alpha})$ is a tensor product of
$W_{1\alpha(1)},\ldots, W_{m\alpha(m)}$ where $\nu_{\alpha}$ is defined by
\[
\nu_{\alpha}(e_{1\gamma(1)}, \ldots, \nu_{m\gamma(m)})\coloneq 
\nu(e_{1\gamma(1)}, \ldots, \nu_{m\gamma(m)}).
\] 
Thus, $P_{\alpha}=\langle \nu_{\alpha}(e_{1\gamma(1)}, \ldots, \nu_{m\gamma(m)}) \Mid \gamma\in \bD_{\alpha}\rangle$. 
Note that 
\[
\dim(P_{\alpha}) = \prod_{i=1}^m \left|D_{i\alpha(i)}\right|=
\prod_{i=1}^m \dim(W_{i\alpha(i)}).
\]
Thus,  by definition~\ref{def:subtenpros}, $(P_{\alpha}, \nu_{\alpha})$ 
is a subspace tensor product of $(P,\nu)$. 
\end{srem}
\begin{theorem}
\label{thm:tnsprddrsms}
Let $V_1, \ldots, V_m$ be vector spaces of dimensions $\dim(V_1)=n_i$.
Suppose $V_i=\oplus_{t=1}^{r_i}W_{it}$, $i=1, \ldots, m$, is the direct sum of subspaces $W_{it}$. 
Then
\begin{equation}
\label{eq:bscsbsntrcng}
\bigotimes_{i=1}^m V_i=\bigotimes_{i=1}^m \bigoplus_{t=1}^{r_i} W_{it} =
\bigoplus_{\alpha\in \Gamma(r_1,\ldots, r_m)}\bigotimes_{i=1}^m W_{i\alpha(i)}\end{equation}
where the $\otimes_{i=1}^m W_{i\alpha(i)}$ are subspace tensor products of 
$\otimes_{i=1}^m V_i.$
\begin{proof}
Let $E_i=\{e_{i1}, \ldots, e_{in_i}\}$, ordered bases of  $V_i$, $i=1, \ldots, m$,  
and $\bD_i=\{D_{i1}, \ldots, D_{ir_i}\}$, ordered partitions of $\underline{n_i}$, be such that
$W_{ij} = \langle e_{it}\Mid t\in D_{ij}\rangle$.
The associated basis for $\otimes_{i=1}^m V_i$ is
$
\{e_{1,\gamma(1)} \otimes \cdots \otimes e_{m,\gamma(m)}\Mid \gamma\in\Gamma(n_1,\ldots, n_m)\}.
$
From definition~\ref{def:prtindnsb}, the induced partition of 
$\Gamma(n_1, \ldots, n_m)$ is 
\[ 
\bD_{\Gamma}=\{\bD_{\alpha}\Mid\bD_{\alpha}=\times_{i=1}^m D_{i\alpha(i)},\; \alpha\in\underline{r_1}\times\cdots\underline{r_m}\}.
\]
Thus, 
\[
\bigotimes_{i=1}^m V_i = 
\bigoplus_{\alpha\in \Gamma(r_1,\ldots, r_m)} \langle e_{1,\gamma(1)} \otimes \cdots \otimes e_{m,\gamma(m)}\Mid 
\gamma\in \bD_{\alpha} \rangle.
\]
Each vector space $P_{\alpha} =\langle e_{1,\gamma(1)} \otimes \cdots \otimes e_{m,\gamma(m)}\Mid \gamma\in \bD_{\alpha} \rangle$ is a subspace of $\otimes_{i=1}^m V_i$ and has dimension
 $\prod_{i=1}^m |D_{i\alpha(i)}| = \prod_{i=1}^m \dim(W_{i\alpha(i)}).$
 If we define $\nu_{\alpha} : \times_{i=1}^m W_{i\alpha(i)} \rightarrow P_{\alpha}$ by
 $\nu_{\alpha}(x_1, \ldots, x_m) \coloneq \nu(x_1, \ldots, x_m)$.  
 Then, $(P_{\alpha}, \nu_{\alpha})$ is the required construction.
\end{proof}
\end{theorem}
%EXAMPLE
\begin{srem}[\bfseries Example: tensor products of direct sums]
\label{srem:exmtnsprddrcsms}
We follow the discussion~\ref{srem:tnsprddrsms}.
Let $(P, \nu)$ be a tensor product of $V_i$, 
$\dim(V_i)=n_i$ with $n_1=n_3=2$ and $n_2=n_4=4.$  Let $E_i=\{e_{i1}, \ldots, e_{in_i}\}$ be ordered bases for the $V_i$, $i=1, \ldots, 4$.
Specifically, take $V_1 =\mathbf{M}_{2,1}$, $V_2 =\mathbf{M}_{2,2}$,
$V_3 =\mathbf{M}_{1,2}$ and $V_4 =\mathbf{M}_{2,2}$. Define the bases
%EQN
\begin{equation}
\label{eq:bsselmv1ex1}
\begin{matrix}
e_{11}=\left(\begin{smallmatrix}1\\0\end{smallmatrix}\right),\;
e_{12}=\left(\begin{smallmatrix}0\\1\end{smallmatrix}\right),\;
e_{31}=\left(\begin{smallmatrix}1&0\end{smallmatrix}\right),\;
e_{32}=\left(\begin{smallmatrix}0&1\end{smallmatrix}\right)
\\
%NEWLINE
e_{21}=e_{41}=
\left(
\begin{smallmatrix}
1\,&0\\
0\,&0
\end{smallmatrix}
\right),\;\;
e_{22}=e_{42}=
\left(
\begin{smallmatrix}
0\,&1\\
0\,&0
\end{smallmatrix}
\right),\;\;
\\
%NEWLINE
e_{23}=e_{43}=
\left(
\begin{smallmatrix}
0\,&0\\
1\,&0
\end{smallmatrix}
\right),\;\;
e_{24}=e_{44}=
\left(
\begin{smallmatrix}
0\,&0\\
0\,&1
\end{smallmatrix}
\right).\;\;
\end{matrix}
\end{equation}
Assume the partitions $D_i$ of $\underline{n_i}$ are
\begin{equation}
D_1=D_3=\{\{1\}, \{2\}\},\; D_2=D_4=\{\{1,2,3,4\}\}. 
\end{equation}
\begin{equation}
V_1=W_{11} \oplus W_{12}\;\mathrm{where}\; W_{11}=\langle e_{11}\rangle, \;W_{12}=\langle e_{12}\rangle\;\mathrm{and}\; V_2=W_{21} = \mathbf{M}_{2,2}
\end{equation}
\begin{equation}
V_3=W_{31} \oplus W_{32}\;\mathrm{where}\; W_{31}=\langle e_{31}\rangle, \;W_{12}=\langle e_{32}\rangle\;\mathrm{and}\; V_4=W_{41} = \mathbf{M}_{2,2}.
\end{equation}
We have
\[ 
\bD_{\Gamma}=\{\bD_{\alpha}\Mid\bD_{\alpha}=\times_{i=1}^4 D_{i\alpha(i)},\; \alpha\in\underline{2}\times\underline{1}\times\underline{2}\times\underline{1}\}.
\]
Thus, $\otimes_{i=1}^4V_i = $
\begin{equation}
\label{srem:exmtnsprddrcsms1}
\begin{matrix}
\langle e_{11}\rangle \otimes \mathbf{M}_{2,2} \otimes \langle e_{31}\rangle
\otimes \mathbf{M}_{2,2} \bigoplus 
\langle e_{11}\rangle \otimes \mathbf{M}_{2,2} \otimes \langle e_{32}\rangle\otimes
\mathbf{M}_{2,2}\bigoplus
\\
\langle e_{12}\rangle \otimes \mathbf{M}_{2,2} \otimes \langle e_{31}\rangle
\otimes \mathbf{M}_{2,2} \bigoplus 
\langle e_{12}\rangle \otimes \mathbf{M}_{2,2} \otimes \langle e_{32}\rangle\otimes
\mathbf{M}_{2,2}.
\end{matrix}
\end{equation}
Each of the summands in equation~\ref{srem:exmtnsprddrcsms1}
is a subspace tensor product (\ref{srem:subtenpros}) of $\otimes_{i=1}^4 V_i,$
and each is isomorphic to $\mathbf{M}_{2,2} \otimes \mathbf{M}_{2,2},$
each with a different isomorphism.
For example, the basis elements of 
$\langle e_{11}\rangle \otimes \mathbf{M}_{2,2} \otimes \langle e_{31}\rangle
\otimes \mathbf{M}_{2,2}$ are
\begin{equation}
\label{eq:bsssmm1}
\{e_{11}\otimes e_{2\beta(2)}\otimes e_{31} \otimes e_{4\beta(4)}\Mid \beta\in \{1,2\}^{\{2,4\}}\}
\end{equation}
corresponding bijectively to basis elements 
$
\{e_{2\beta(2)} \otimes e_{4\beta(4)}\Mid \beta\in \{1,2\}^{\{2,4\}}\}
$
of $\mathbf{M}_{2,2} \otimes \mathbf{M}_{2,2}.$
\end{srem}
\begin{srem}[\bfseries Example: tensor products of direct sums -- Kronecker products]
\label{srem:exmtnsprddrcsmskrn}
In example~\ref{srem:exmtnsprddrcsms} our model for the tensor products defined the underlying vector spaces as matrices but didn't use any properties of them except for the dimensions.
Let $(P, \nu)$ be a tensor product of $V_i$, $i=1, \ldots, 4$,
$\dim(V_i)=n_i$ with $n_1=n_3=2$ and $n_2=n_4=4.$  
Specifically, take $V_i\in \mathbf{M}_{p_i,q_i}$ as follows: $V_1 =\mathbf{M}_{2,1}$, $V_2 =\mathbf{M}_{2,2}$, $V_3 =\mathbf{M}_{1,2}$ and $V_4 =\mathbf{M}_{2,2}$. 
Take $P=\mathbf{M}_{p,q}$ where $p=p_1p_2p_3p_4=8$ and $q=p_1q_2q_3q_4=8$.
In this example, we associate the bases of these vector spaces with the standard bases of matrices (as used in discussions~\ref{rem:tnsprdmtr}, \ref{rem:unvfctprpexm}, \ref{rem:tnsprdmtr2}, and \ref{rem:krnprd}).
\begin{equation}
\label{eq:bsselmvkrnprd1}
\begin{matrix}
E^{(1)}_{11}=\left(\begin{smallmatrix}1\\0\end{smallmatrix}\right),\;
E^{(1)}_{21}=\left(\begin{smallmatrix}0\\1\end{smallmatrix}\right),\;
E^{(3)}_{11}=\left(\begin{smallmatrix}1&0\end{smallmatrix}\right),\;
E^{(3)}_{12}=\left(\begin{smallmatrix}0&1\end{smallmatrix}\right)
\\
%NEWLINE
E^{(2)}_{11}=E^{(4)}_{11}=
\left(
\begin{smallmatrix}
1\,&0\\
0\,&0
\end{smallmatrix}
\right),\;\;
E^{(2)}_{12}=E^{(4)}_{12}=
\left(
\begin{smallmatrix}
0\,&1\\
0\,&0
\end{smallmatrix}
\right),\;\;
\\
%NEWLINE
E^{(2)}_{21}=E^{(4)}_{21}=
\left(
\begin{smallmatrix}
0\,&0\\
1\,&0
\end{smallmatrix}
\right),\;\;
E^{(2)}_{22}=E^{(4)}_{22}=
\left(
\begin{smallmatrix}
0\,&0\\
0\,&1
\end{smallmatrix}
\right).\;\;
\end{matrix}
\end{equation}
Each $E^{(i)}$ is ordered lexicographically.  
Thus, $E^{(2)}=(E^{(2)}_{11}, E^{(2)}_{12}, E^{(2)}_{21}, E^{(2)}_{22})$.
In this case, equation~\ref{srem:exmtnsprddrcsms1} becomes $\otimes_{i=1}^4V_i = $
\begin{equation}
\label{eq:bsselmvkrnprdkpr2}
\begin{matrix}
\langle E^{(1)}_{11}\rangle \otimes \mathbf{M}_{2,2} \otimes \langle E^{(3)}_{11}\rangle
\otimes \mathbf{M}_{2,2} \bigoplus 
\langle E^{(1)}_{11}\rangle \otimes \mathbf{M}_{2,2} \otimes \langle E^{(3)}_{12}\rangle\otimes
\mathbf{M}_{2,2}\bigoplus
\\
\langle E^{(1)}_{21}\rangle \otimes \mathbf{M}_{2,2} \otimes \langle E^{(3)}_{11}\rangle
\otimes \mathbf{M}_{2,2} \bigoplus 
\langle E^{(1)}_{21}\rangle \otimes \mathbf{M}_{2,2} \otimes \langle E^{(3)}_{12}\rangle\otimes
\mathbf{M}_{2,2}.
\end{matrix}
\end{equation}
We interpret these tensor products as Kronecker products of matrices as in 
examples~\ref{rem:krnprd} and~\ref{srem:krnprdexm}.
Let $A_i=(a_i(s,t))\in \mathbf{M}_{p_i,q_i}$, $i=1, \dots, 4$.
We have 
\[
A_i =
\sum_{(\mu(i), \kappa(i))\in \Gamma (p_i,q_i)}a_i(\mu(i), \kappa(i)) E^{(i)}_{\mu(i), \kappa(i)}
\;\mathrm{and}
\] 
\begin{equation}
\label{eq:krnprd4}
A_1\otimes A_2 \otimes A_3\otimes A_4 =\sum_{\mu,\kappa} \prod_{i=1}^4 a_i(\mu(i),\kappa(i)) E_{\mu,\kappa}
\end{equation}
where the $E_{\mu,\kappa}(\alpha, \beta)=\mathcal{X}((\alpha,\beta)=(\mu,\kappa))$
are the standard basis elements of $\mathbf{M}_{p,q}$. 
If $\mathbf{A} = A_1\otimes \cdots \otimes A_4$, the entry 
$\mathbf{A}(\mu, \kappa) = \prod_{i=1}^4 a_i(\mu(i),\kappa(i))$.
Denote the sequences of possible row values by  $\Gamma_{\mathbf{r}} \coloneq \Gamma(p_1,p_2,p_3,p_4)=\Gamma(2,2,1,2)$ and  the sequences of column values by
$\Gamma_{\mathbf{c}}\coloneq \Gamma(q_1,q_2,q_3,q_4)=\Gamma(1,2,2,2)$. 
Order both $\Gamma_{\mathbf{r}}$ and $\Gamma_{\mathbf{c}}$ lexicographically.
The basis matrices
\begin{equation}
\label{eq:krnprd5}
\{E_{\mu,\kappa} \Mid \mu\in \Gamma_{\mathbf{r}}\,,\,  \kappa\in \Gamma_{\mathbf{c}}\}
\end{equation}
are ordered 
by lexicographic order on $(\mu, \kappa) \in \Gamma_{\mathbf{r}}\times \Gamma_{\mathbf{c}}.$
In matrix~\ref{eq:rwsclmslex}, the rows are shown indexed lexicographically by 
$\Gamma_{\mathbf{r}}$ and the columns by $\Gamma_{\mathbf{c}}$
(corresponding to $\Phi=\Gamma_{\mathbf{r}}$ and $\Lambda=\Gamma_{\mathbf{c}}$ in notation~\ref{srem:addcon} below ).  
Note that if 
\[
A_1\otimes \cdots \otimes A_4=\mathbf{A} \in 
\langle E^{(1)}_{11}\rangle \otimes \mathbf{M}_{2,2} \otimes \langle E^{(3)}_{11}\rangle
\otimes \mathbf{M}_{2,2}
\] 
then the support of $A$ (\ref{def:mat}) is in the submatrix of~\ref{eq:rwsclmslex} labeled with ``$a$''.
Similarly for 
\[
B_1\otimes \cdots \otimes B_4=\mathbf{B} \in 
\langle E^{(1)}_{11}\rangle \otimes \mathbf{M}_{2,2} \otimes \langle E^{(3)}_{12}\rangle
\otimes \mathbf{M}_{2,2}\; (\mathrm{labeled}\;b),
\] 
\[
C_1\otimes \cdots \otimes C_4=\mathbf{C} \in 
\langle E^{(1)}_{21}\rangle \otimes \mathbf{M}_{2,2} \otimes \langle E^{(3)}_{11}\rangle
\otimes \mathbf{M}_{2,2}\; (\mathrm{labeled}\;c),
\]
\[
D_1\otimes \cdots \otimes D_4=\mathbf{D} \in 
\langle E^{(1)}_{21}\rangle \otimes \mathbf{M}_{2,2} \otimes \langle E^{(3)}_{12}\rangle
\otimes \mathbf{M}_{2,2}\; (\mathrm{labeled}\;d).
\]
\begin{equation}
\label{eq:rwsclmslex} 
\bordermatrix{~ &1111 & 1112 &1121&1122&1211 & 1212 &1221&1222 \cr
%L1
1111& a & a &b & b& a & a &b& b  \cr
%L2
1112&  a & a &b & b& a & a &b & b  \cr 
%L3
1211&  a & a &b & b & a & a &b & b \cr
%L4 
1212&  a & a &b & b & a & a &b & b \cr
%L1
2111& c & c &d& d & c & c &d & d \cr
%L2
2112&  c & c &d & d & c & c &d & d \cr 
%L3
2211&  c & c &d & d & c & c &d & d \cr
%L4 
2212&  c & c &d & d & c & c &d & d \cr
}
\end{equation}
\end{srem}

%SECTION    
\section{Background concepts and notation}
\label{sec:review}
In this section we review concepts and notation from basic discrete mathematics courses.  This section can be skipped and reviewed as needed.  Wikipedia is a good source.
\begin{srem}[\bfseries Sets and lists]
\label{rem:notation}
The empty set is denoted by $\emptyset$. Sets are specified by braces: $A=\{1\}$, $B=\{1,2\}$.  They are unordered, so
$B=\{1,2\}=\{2,1\}$.  Sets $C=D$ if  $x\in C\; \mathrm{implies}\; x\in D$ 
(equivalently, $C\subseteq D$) and
$x\in D\; \mathrm{implies}\; x\in C$.  If you write $C=\{1,1,2\}$ and $D=\{1,2\}$ then, by the definition of set equality, $C=D$.  
If $A$ is a set then $\bP(A)$ is all subsets of $A$ and $\bP_k(A)$ is all subsets of cardinality, $|A|=k$.

A {\em list}, {\em vector}  or {\em sequence} (specified by parentheses) is ordered:
$C' = (1,1,2)$ is not the same as $(1,2,1)$ or $(1,2)$.  Two lists (vectors, sequences), 
$(x_1, x_2, \ldots , x_n) = (y_1, y_2, \ldots , y_m)$, are equal if and only if $n=m$ and
$x_i = y_i$ for $i=1, \ldots , n$.  
A list such as $(x_1, x_2, \ldots , x_n)$ is also written $x_1, x_2, \ldots , x_n$, without the parentheses.  Sometimes a list $L$ will be specified as a set,
$S$ with additional information defining the linear 
order $<_S.$ Lex order~\ref{def:lxcrdr} is an example.
\end{srem}
 %DEF
\begin{definition}[\bfseries Lexicographic, ``lex,'' order]
\label{def:lxcrdr}
Let $C_i$, $i=1,\ldots, k$, be lists of distinct elements, and let
$L=C_1\times \cdots \times C_k$ be the their product as sets.
Define lexicographic order on $L$, indicated by $<_L$,  by
$(a_1, \ldots, a_k) <_L (b_1, \ldots, b_k)$ if $a_1< b_1$ or if there is some 
$t\leq k$ such that $a_i = b_i$, $i<t$, but $a_t < b_t$.
We have used ``$<$'' for the various linear orders on the $C_i$ and will usually do the same for lexicographic (``lex'') order, replacing $<_L$ by simply $<$.
\end{definition}
%DEF
\begin{definition}[\bfseries Position and rank functions for linear orders]
\label{rem:pstrnk}
Let $(\Lambda, \leq)$ be a finite linearly ordered set.
For $x\in \Lambda$ and $S\subseteq \Lambda$ let
\begin{equation}
\pi^{\Lambda}_{S}(x)\coloneq\left|\{t\Mid t\in S, t\leq x\}\right|\;\;\mathrm{and}\;\;
\rho^{\Lambda}_{S}(x)\coloneq\left|\{t\Mid t\in S, t < x\}\right|.
\end{equation}
$\pi^{\Lambda}_{S}$ is called the {\em position function} for $\Lambda$ relative  to $S$ and 
$\rho^{\Lambda}_S$ the {\em rank function} for $\Lambda$ relative to $S$.
If $S=\Lambda$, we use $\pi^{\Lambda}$ instead of $\pi^{\Lambda}_{\Lambda}$, and,
similarly, we use $\rho^{\Lambda}$ instead of $\rho^{\Lambda}_{\Lambda}$.
\end{definition}
%DEF
\begin{definition} [\bfseries Order relation, partially ordered set, poset]
\index{order relation}\index{relation, order}
\label{def:ordrel}
A subset $R\subseteq S\times S$ is called a {\em binary relation on} $S$.
The statement $(x,y)\in R$ is also denoted by $x\mathrel R y$.
Likewise, $(x,y)\notin R$ is denoted by $x\not\mathrel R y$.
A binary relation on a set $S$ is called an {\em order relation\/} if
it satisfies the following three conditions and (usually
written $x\preceq y$ instead of $x\mathrel R y$ in this case).
\begin{itemize}
\item {(i)} (Reflexive)
\index{reflexive relation}\index{relation, reflexive}
For all $s\in S$ we have $s\preceq s$.
\item {(ii)} (Antisymmetric)
\index{antisymmetric relation}\index{relation, antisymmetric}
For all $s,t\in S$ such that $s\neq t$, if $s\preceq t$ then
$t\not\preceq s$.
\item {(iii)} (Transitive)
\index{transitive relation}\index{relation, transitive}
For all $r,s,t\in S$,  $r\preceq s$ and $s\preceq t$ implies $r\preceq t$.
\end{itemize}

\index{set, partially ordered}
\index{poset}
\noindent
A set $S$ together with an order relation $\preceq$ is called a
{\em partially ordered set\/} or {\em poset}, $(S, \preceq)$.
If for all $(x,y)\in S\times S$, either $x\preceq y$ or $y\preceq x$ (both if $x=y$) then $(S, \preceq)$ is a {\em linearly ordered set} 
(e.g., lists~\ref{rem:notation}, lexicographic order~\ref{def:lxcrdr}).
\end{definition}
\begin{defn}[\bfseries \bf Function] 
\label{def:function}
Let $A$ and $B$ be sets.   
A function $f$ from $A$ to $B$ is a rule that assigns to each element 
$x\in A$ a unique element $y \in B$.  We write $y=f(x)$.
Two functions $f$ and $g$ from $A$ to $B$ are equal if $f(x)=g(x)$ for all $x\in A$.
\end{defn}
Given a function $f$ from $A$ to $B$, we can define a set $F\subseteq A\times B$ by \index{function!graph}
\begin{equation}
\label{eq:graphfunction}
F=\{(x,f(x))\Mid x\in A\}
\end{equation}
We call $F$ the {\em graph} of $f$, denoted by ${\rm Graph}(f)$. 
 A subset $F\subseteq A\times B$ is the graph  of a function from $A$ to $B$ if and only if it satisfies 
 the following two conditions:
\begin{equation}
{\rm G1}:\;\;(x,y)\in F\;\;\;  {\rm and}\;\;\; (x,y')\in F \implies y=y' 
\end{equation}
\begin{equation}
{\rm G2}:\;\;\{x\Mid (x,y)\in F\} = A.
\end{equation}
Two functions, $f$ and $g$, are equal if and only if their graphs are equal as sets:
${\rm Graph}(f) = {\rm Graph}(g).$
\index{function!domain, range}
The set $A$ is called the {\em domain} of $f$ (written $A={\rm domain}(f)$), and $B$ is called the {\em range} of $f$ (written $B={\rm range}(f)$). 
The notation $f:A \rightarrow B$ is used to denote that $f$ is a function with domain $A$
and range $B$.
\index{function! image}
For $S\subseteq A$, define $f(S)$ (image of $S$ under $f$) by
$f(S)\equiv \{f(x)\Mid x \in S\}$. 
In particular, $f(A)$ is called the {\em image} of $f$ (written $f(A)={\rm image}(f)$). 
The set of {\em all} functions with domain $A$ and range $B$ can be written  $\{f\Mid f:A \rightarrow B\}$ or simply as $B^A$.  If $A$ and $B$ are finite then $\left| B^A\right|$ is  $|B|^{|A|}$.
\index{function!indicator or characteristic}
The {\em characteristic or indicator function} of a set $S\subseteq A$, 
$\mathcal{X}_S:A\rightarrow \{0,1\}$, is defined by
\begin{equation}
\label{eq:charfunc}
\mathcal{X}_S(x) = 1 \,\,{\rm if\;and\;only\;if\,\,} x\in S.
\end{equation} 
\index{function!restriction, composition}
The {\em restriction} $f_S$ of  $f:A\rightarrow B$ to a subset $S\subseteq A$ is defined by
\begin{equation}
\label{eq:restriction}
f_S: S \rightarrow B\;\;{\rm where}\;\;f_S(x)=f(x)\;\;{\rm for\;\;all\;\;} x\in S.
\end{equation}

If $f:A\rightarrow B$ and $g:B \rightarrow C$ then the {\em composition} of $g$ and $f$, denoted
by $gf: A\rightarrow C$, is defined by 
\begin{equation}
\label{eq:composition}
gf(x) = g(f(x))\;\; {\rm for}\;\; x \in A.
\end{equation}
There are many ways to describe a function.  Any such description must specify the domain, the range, and the rule for assigning some range element to each domain element.
You could specify a function using set notation:  the domain is  $\{1,2,3,4\}$, the range is $\{a, b, c, d, e\}$, and the function $f$ is the set: 
${\rm Graph}(f)=\{(1,a)\,(2,c)\,(3,a)\,(4,d)\}$. 
\index{function!two line description} 
Alternatively, you could describe the same function by giving the range as $\{a, b, c, d, e\}$ and using {\em two line notation}
\begin{equation}
\label{eq:twolinedis}
f=
\left(
\begin {array}{cccc}
1 & 2 & 3 & 4 \\
a & c & a & d 
\end{array}
 \right).
\end{equation}
If we assume the domain, in order, is $1\,2\,3\,4$, then~\ref{eq:twolinedis} can be abbreviated  to {\em one line}: $a\,c\,a\,d$.
Sometimes it is convenient to describe a function $f\in B^A$ by the ``maps to'' notation.  For example, $2\mapsto c$, $4\mapsto d$, otherwise $x\mapsto a$.
%SETS OF FUNCTIONS
\index{functions!strictly increasing SNC}
\index{functions!weakly increasing WNC}
\index{functions!injective SNC}
\index{functions!permutations PER}
\begin{defn}[\bfseries Sets of functions] 
\label{def:setsfunc}
Let ${\underline n}=\{1, 2, \ldots, n\}$ and let
${\underline p}^{\underline n}$ be all functions with domain  
${\underline n}$, range ${\underline p}$.  Define
\[ 
{\rm SNC}(n,p) = \{f\Mid f\in {\underline p}^{\underline n}, i<j \implies f(i)<f(j)\}
\;\;{\bf (strictly\;increasing)}
\]
\[ 
{\rm WNC}(n,p) = \{f\Mid f\in {\underline p}^{\underline n}, i<j \implies f(i)\leq f(j)
\;\;{\bf (weakly\;increasing)}
\]
\[ 
{\rm INJ}(n,p) = \{f\Mid f\in {\underline p}^{\underline n}, i\neq j \implies f(i)\neq f(j)\}
\;\;{\bf (injective)}
\]
\[
{\rm PER}(n) = {\rm INJ}(n,n) 
\;\;{\bf (permutations\;of\;}{\underline n}).
\] 
From combinatorics, $|{\rm INJ}(n,p)|=(p)_n=p(p-1)\cdots (p-n+1)$, |PER(n)| = n!,
\[
|{\rm SNC}(n,p)|= \left(\begin{array}{c} p\\n\end{array}\right)
\;\;{\rm and}\;\;
|{\rm WNC}(n,p)|= \left(\begin{array}{c} p+n-1\\n\end{array}\right).    
\]
More generally, if $X\subseteq {\underline n}$ and $Y\subseteq {\underline p}$, then
${\rm SNC}(X,Y)$ denotes the strictly increasing functions from $X$ to $Y$.  We define
${\rm WNC}(X,Y)$ and ${\rm INJ}(X,Y)$ similarly.
Sometimes  ``increasing'' is used instead of  ``strictly increasing'' or
``nondecreasing''   instead of  ``weakly increasing''
\end{defn}
%DEF 
\begin{definition}[\bfseries Matrix]
\label{def:mat}
Let $m, n$ be positive integers.  An $m$ by $n$ matrix with entries in a set $S$ is a function 
$f: \underline{m}\times \underline{n} \rightarrow S$.  The sets $\underline{m}$ and $\underline{n}$ are the {\em row indices} and {\em column indices} respectively.
The set of all such $f$ is denoted by $\mathbf{M}_{m,n}(S)$.
If $S=F$, a field, the set $\{x\mid f(x)\neq 0\}$ is the {\em support} of $f$, and the matrices 
$E_{ij}\in\mathbf{M}_{m,n}(F)$ with a $1$ in position $i,j$ and $0$ elsewhere are the {\em standard basis elements}.  Instead of $E_{ij}$, we use $E_{i,j}$ (with comma) when needed for clarity.
\end{definition}
%SREM
\begin{srem}[\bfseries Additional matrix notational conventions]
\label{srem:addcon}
More generally, a matrix is a function $f:\Phi\,\times\, \Lambda \rightarrow S$ where $\Phi$ and $\Lambda$ are linearly ordered sets (row and column indices respectively).
We use $A[X\vert Y]$ to denote the submatrix of $A$ gotten by retaining rows indexed by the set $X\subseteq \Phi$ and columns indexed by the set 
$Y\subseteq \Lambda$.
We use $A(X\vert Y)$ to denote the submatrix of $A$ gotten by retaining rows indexed by the set $\Phi\setminus X$ (the complement of $X$ in $\Phi$) and columns indexed by the set $\Lambda \setminus Y$.  
 We also use the mixed notation  $A[X\vert Y)$ and  $A(X\vert Y]$
 with obvious meaning.
 We use $\Theta$ to denote the zero matrix of the appropriate size and $I$ to denote the identity matrix.
\end{srem}

%DEF
\begin{defn}[\bfseries Partition of a set]
\label{def:setpartition}
A partition of a set $Q$ is a collection, ${\mathcal B}(Q)$, of nonempty subsets, 
$X$, of $Q$ such that each element of $Q$ is contained in {\em exactly one}  set
$X\in{\mathcal B}(Q).$ 
The sets $X\in{\mathcal B}(Q)$ are called the {\em blocks} of the partition
 ${\mathcal B}(Q).$   
A set $D\subseteq Q$ consisting of exactly one   element from each block is called a  {\em system of distinct representatives} (or ``SDR'') for the partition.
If each $X\in \mathcal {B}(Q)$ has $|X|=1$, then we call ${\mathcal B}(Q)$ the
{\em discrete} partition.  If the partition $\mathcal {B}(Q)=\{Q\}$ is called the
{\em unit} partition.  The set of all partitions of $Q$ is $\Pi(Q)$.
The set of all partitions of $Q$ with $k$ blocks is $\Pi_k(Q)$.
\end{defn}
If $|Q|=n,$ the numbers $S(n,k)\coloneq |\Pi_k(Q)|$ are called the {\em Stirling numbers of the second kind}  and the numbers, $B(n)\coloneq |\Pi(Q)|$   are called the Bell numbers.
If $Q=\{1,2,3,4,5\}$ then $\mathcal{B}(Q)=\{\{1,3,5\},\{2,4\}\}\in \Pi_2(Q)$ is a partition of $Q$ with two blocks: $\{1,3,5\}\in \mathcal{B}(Q)$ and $\{2,4\}\in \mathcal{B}(Q).$
The set $D=\{1,2\}$ is an SDR for $\mathcal{B}(Q)$.
$S(5,2)=15$ and $B(5)=52$.    
\index{function!coimage}
%DEF
\begin{defn}[\bfseries Coimage partition]
\label{def:coimage}
Let $f:A\rightarrow B$ be a function with domain $A$ and range $B$.  
Let ${\rm image}(f) = \{f(x)\Mid x\in A\}$ (\ref{def:function}).  The inverse image of an
element $y\in B$ is the set $f^{-1}(y)\equiv \{x\Mid f(x)=y\}.$
The {\em coimage} of  $f$ is the set of subsets of $A:$ 
\begin{equation}
{\rm coimage}(f) = \{f^{-1}(y)\Mid y\in {\rm image}(f)\}.
\end{equation}
The ${\rm coimage}(f)$ is a partition of $A$ (\ref{def:setpartition}) called the
{\em coimage partition of $A$ induced by }$f$.
\end{defn}

For the function $f$ of \ref{eq:twolinedis}, we have ${\rm image}(f) = \{a,c,d\}.$
Thus, the coimage of $f$ is 
\begin{equation}
\label{eq:examplecoimage}
{\rm coimage}(f) = \{f^{-1}(a), f^{-1}(c), f^{-1}(d)\}=\{\{1,3\},\{2\},\{4\}\}.
\end{equation}
\hspace*{1 in}

%DEF
\begin{srem}[\bfseries Posets of subsets and partitions]
\label{srem:pstsbsprt}
If $S=\bP(Q)$ denotes all subsets of $Q$, then $(S, \preceq)$ is a poset
(\ref{def:ordrel}) if $\preceq \coloneq \subseteq$ (set inclusion).  In this case, $(S, \subseteq)$ is called the {\em poset of subsets} or {\em lattice of subsets} of $Q$. 

Suppose $S=\Pi(Q)$ is the collection of all partitions, ${\mathcal B}(Q)$, of $Q$.  
Let $\preceq$ denote {\em refinement} of partitions where  $\mathcal{B}_1(Q)\preceq \mathcal{B}_2(Q)$ means the blocks of 
$\mathcal{B}_1(Q)$ are obtained by further subdividing the blocks of 
$\mathcal{B}_2(Q)$ (e.g., if $\mathcal{B}_1(Q)=\{\{1,3\},\{5\},\{2\},\{4\}\}$ and 
$\mathcal{B}_2(Q)=\{\{1,3,5\},\{2,4\}\}$  then $\mathcal{B}_1(Q)\preceq \mathcal{B}_2(Q)$).
Using general notation for posets, if  $x=\mathcal{B}_1$ and $y=\mathcal{B}_2$, $x,y\in {\mathcal B}(Q)$, 
$x\prec y$, then 
$x$ {\em is covered by} $y$, written 
 $x\prec_c y$, if $\{z\Mid x\prec z \prec y\}=\emptyset$.  
 Diagram~\ref{srem:Fig003Hasse4} represents the poset 
$(\Pi(\underline{4}),\preceq)$ by its {\em covering relation} where
  $x\prec_c y$ is represented by $y\rightarrow x$.
Such a diagram for a poset is called a {\em Hasse diagram}.

 \end{srem}
%MINIPAGE HASSE4
\noindent
\begin{minipage}{\textwidth}%
\begin{srem}[\bfseries Diagram of covering relation for $\Pi(\underline{4})$]
\label{srem:Fig003Hasse4}
\end{srem}
%\begin{equation}
%\mathbf{Diagram\;of\;covering\;relation\;for\;\Pi(\underline{4})}
%\end{equation}
\begin{center}
\includegraphics[scale=.8]{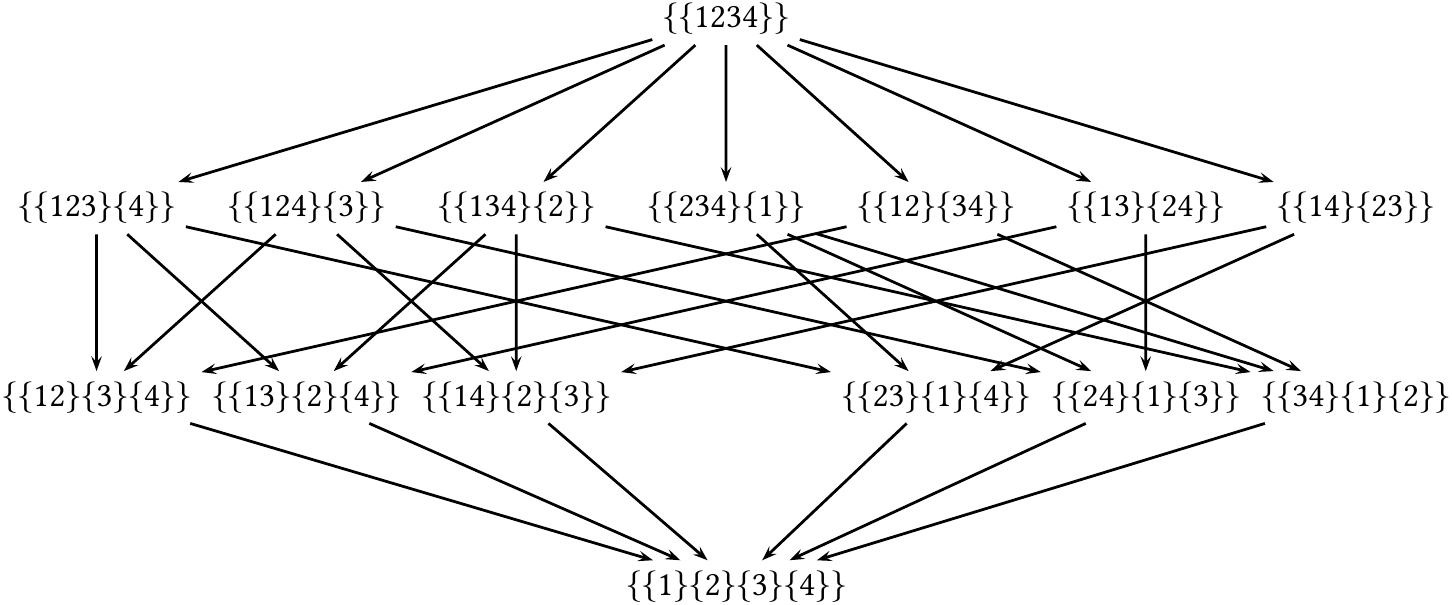}
\end{center}
\end{minipage}
%END MINIPAGE

\section{Discussion and acknowledgements}
\label{sec:discussion}
This article is a rewrite of notes I originally prepared for first year graduate students in combinatorics seminars in which multilinear algebra was applied to combinatorics (UCSD Mathematics and CSE).   Additional references in multilinear algebra for these seminars were provided by Professor Marvin Marcus, my thesis advisor, friend and mentor.

S. Gill Williamson, 2015\\
\url{http://cseweb.ucsd.edu/~gill}
\bibliographystyle{alpha}
\bibliography{Tensors}

%S. Gill Williamson, 2012\\
%${\rm cseweb.ucsd.edu\slash\sim gill\slash}$  
%
%\tableofcontents
%\index{contents}
%\hyperlink{index}{Index}
%\newpage
%\thispagestyle{empty}

%\label{cor:relpridiaent}
%\newpage
%\mbox{}
%\newpage
%\hypertarget{index}{ }
%\printindex
%\newpage
%\centerline{NOTES}
%\newpage
%\centerline{NOTES}
%\newpage
%\centerline{NOTES}
\end{document}